\documentclass[11pt]{article}
\title{The Local Lemma is asymptotically tight for SAT\thanks{A preliminary version
    of this paper has appeared as an extended abstract in the
    Proceedings of SODA 2011}}
\author{H. Gebauer
\thanks{Institute of Applied Mathematics and Physics, Zurich University of 
Applied Sciences, Switzerland; mail: geba@zhaw.ch.},
 T. Szab\'o
 \thanks{Department of Mathematics and Computer Science, Freie Universit\"{a}t Berlin, 14195 Berlin, Germany; mail: szabo@math.fu-berlin.de.
 Research supported by DFG, project SZ 261/1-1.},
 G. Tardos \thanks{R\'{e}nyi
   Institute, Budapest, Hungary; mail: tardos@renyi.hu. Supported by the
   Cryptography Lend\"ulet project of the Hungarian Academy of Sciences and the
 Hungarian grant OTKA K-116769}
 }
 \date{}

\usepackage{amsmath, amsfonts}
\usepackage{amsbsy}
\usepackage{amsgen}
\usepackage{amssymb}
\usepackage{amsxtra}
\usepackage{amsthm}
\usepackage{graphicx}

\newtheorem{thm}{Theorem}[section]
\newtheorem{cor}[thm]{Corollary}

\newtheorem{problem}[thm]{Problem}
\newtheorem{lemma}[thm]{Lemma}

\newtheorem{obse}[thm]{Observation}
\newtheorem{prop}[thm]{Proposition}
\theoremstyle{definition}

\theoremstyle{remark}

\numberwithin{equation}{section}

\newcommand{\Real}{\mathbb R}

\newcommand{\N}{\mathbb N}

\newcommand{\pairing}{\mathsf{pairing}}
\newcommand{\tree}{\mathsf {tree}}

\def\f{f_{\rm tree}}
\usepackage{enumerate}

\begin{document}
\maketitle
\thispagestyle{empty}

\begin{abstract}
The Local Lemma is a fundamental tool of probabilistic combinatorics and theoretical computer science, yet
there are hardly any natural problems known where it provides an asymptotically tight answer. The main theme
of our paper is to identify several of these problems, among them a couple of
widely studied extremal functions related to certain restricted versions of the $k$-SAT problem,
where the Local Lemma does give essentially optimal answers.

As our main contribution, we construct unsatisfiable $k$-CNF formulas where every clause has $k$ distinct literals
and every variable appears in at most $\left( \frac{2}{e} + o(1)\right)\frac{2^k}{k}$ clauses.
The Lopsided Local Lemma, applied with an assignment of random values
according to counterintuitive probabilities,
shows that this is asymptotically best possible.
The determination of this extremal function
is particularly important as it represents the value
where the corresponding $k$-SAT problem exhibits a
complexity hardness jump: from having every instance being a
YES-instance it becomes NP-hard just by
allowing each variable to occur in one more clause.

The construction of our unsatisfiable CNF-formulas is based on the binary
tree approach of \cite{Gebauer} and thus the constructed formulas are in the
class MU(1) of minimal unsatisfiable formulas having one more clauses than
variables. The main novelty of our approach here comes in setting up an appropriate
continuous approximation of the problem. This leads us to a differential equation, the
solution of which we are able to estimate.
The asymptotically optimal binary trees are then obtained through a discretization of
this solution.

The importance of the binary trees constructed is also underlined by their appearance
in many other scenarios. In particular, they give asymptotically precise
answers for seemingly unrelated problems like the European Tenure Game
introduced by Doerr \cite{Doerr04} and a search problem allowing  a
limited number of consecutive lies. 
As yet another consequence we slightly improve the best known bounds
on the maximum degree and maximum edge-degree of a $k$-uniform
Maker's win hypergraph in the Neighborhood Conjecture of Beck.

\end{abstract}

\maketitle
\setcounter{page}1

\section{Introduction}

The satisfiability of Boolean formulas is the archetypal NP-hard
problem.
Somewhat unusually we define a $k$-CNF formula as the conjunction of clauses
that are the disjunction of {\em exactly} $k$ distinct literals. (Note that most
texts allow shorter clauses in a $k$-CNF formula, but fixing the exact length
will be important for us later on.)
The problem of deciding whether a $k$-CNF formula is satisfiable is denoted by
$k$-SAT, it is solvable in polynomial time for $k=2$, and is NP-complete for
every $k\ge3$ as shown by Cook \cite{Cook}.

Papadimitriou and Yannakakis \cite{PY} have shown that MAX-$k$-SAT (finding the
maximum number of simultaneously satisfiable clauses in an input $k$-CNF
formula) is even MAX-SNP-complete for every $k\geq 2$.

The first level of difficulty in satisfying a CNF formula arises
when two clauses share variables.
For a finer view into the transition to NP-hardness, a grading of
the class of $k$-CNF formulas can be introduced, that limits how much
clauses interact locally.
A $k$-CNF formula is called a $(k,s)$-CNF formula if every variable appears in
at most $s$ clauses.
The problem of deciding satisfiability of a $(k,s)$-CNF formula is denoted by
$(k,s)$-SAT, while finding the maximum number of simultaneously satisfiable
clauses in such a formula is called MAX-$(k,s)$-SAT.

Tovey \cite{Tovey} proved that while every $(3,3)$-CNF formula is satisfiable (due to Hall's
theorem), the problem of deciding whether a $(3,4)$-CNF formula is satisfiable
is already NP-hard. Dubois \cite{Dubois90} showed that $(4,6)$-SAT and
$(5,11)$-SAT are also NP-complete.

Kratochv\'{i}l, Savick\'y, and Tuza \cite{KST} defined the value $f(k)$ to be
the largest integer $s$ such that every $(k,s)$-CNF is satisfiable.
They also generalized Tovey's  result by showing that for every
$k \geq 3$ $(k, f(k) + 1)$-SAT is already NP-complete.
In other words, for every $k\geq 3$ the $(k, s)$-SAT problem goes through a
kind of ``complexity phase transition'' at the value $s=f(k)$. On the one hand
the $(k,f(k))$-SAT problem is trivial
by definition in the sense that every instance of the problem is a
``YES''-instance. On the other hand the $(k, f(k)+1)$-SAT problem is already
NP-hard, so the problem becomes hard from being trivial just by allowing one
more occurrence of each variable. For large values of $k$ this might seem
astonishing, as the value of the transition is exponential in $k$: one might
think that the change of just one in the parameter should have hardly any
effect.

The complexity hardness jump is even greater:
MAX-$(k, s)$-SAT is also MAX-SNP-complete for every $s>f(k)$, $k\ge2$
as was shown by Berman, Karpinski, and Scott \cite{BK, BKS} (generalizing a
result of Feige \cite{Feige98} who showed that MAX-$(3,5)$-SAT is hard to
approximate within a certain constant factor).

The determination of where this complexity hardness jump occurs is the topic
of the current paper.

For a lower bound the best tool available is the Lov\'asz Local Lemma. The lemma does
not deal directly with number of occurrences of variables, but rather with
pairs of clauses that share at least one variable. We call such a pair an
{\em intersecting pair} of clauses. A
straightforward consequence of the lemma states that if every clause of a
$k$-CNF formula intersects at most $2^k/e-1$ other clauses, then the formula
is satisfiable. It is natural to ask how tight this bound is and for that
Gebauer et al.\ \cite{GMSW} define $l(k)$ to be the largest integer number
satisfying that whenever all clauses of a $k$-CNF formula intersect at most
$l(k)$ other clauses the formula is satisfiable. With this notation the
Lov\'asz Local Lemma implies that
\begin{equation}
l(k)\ge\left\lfloor\frac{2^k}e\right\rfloor-1.\label{lll}
\end{equation}

The order of magnitude of this bound is trivially optimal: $l(k)<2^k-1$ follows
from the unsatisfiable $k$-CNF formula consisting of all possible $k$-clauses on only
$k$ variables.

In \cite{GMSW} a hardness jump is proved for the function $l$: for $k\ge3$,
deciding the satisfiability of $k$-CNF formulas with maximum neighborhood size
at most $l(k)+2$ is NP-complete.\footnote{In \cite{GMSW} the
  slightly more complicated formula $\max\{ l(k)+2, k+3\}$ appeared for the
  maximum neighborhood size but this can be simplified to just its first term.
For $k\geq 5$ this was already observed in
  \cite{GMSW} to follow
  from Equation~\ref{lll}.
The statement follows from a slightly stronger form of the
  Local Lemma for $k=4$ and from a case analysis for $k=3$.}

As observed by Kratochv\'\i l, Savick\'y and Tuza \cite{KST} the bound
(\ref{lll}) immediately implies
\begin{equation}
f(k)\ge\left\lfloor\frac{l(k)}k\right\rfloor+1\ge\left\lfloor\frac{2^k}{ek}\right\rfloor.\label{kst}
\end{equation}

From the other side Savick\'y and Sgall \cite{SS} showed that
$f(k) = O\left(k^{0.74} \cdot \frac{2^{k}}{k}\right)$.
This was improved by Hoory and Szeider \cite{HSfa} who came within a logarithmic factor:
$f(k) = O\left(\log k \cdot \frac{2^{k}}{k}\right)$.
Recently, Gebauer \cite{Gebauer} showed that the order of
magnitude of the lower bound is correct and $f(k) = \Theta(\frac{2^k}{k})$.

More precisely, the construction of \cite{Gebauer} gave
$f(k) \leq \frac{63}{64}\cdot \frac{2^k}{k}$ for infinitely many $k$.
The constant factor $\frac{63}{64}$ was clearly not the optimum rather
the technical limit of the approach of \cite{Gebauer}.
Determining $f(k)$ asymptotically remained an outstanding open problem and
there was no clear consensus about where the correct asymptotics
should fall between the constants $1/e$ of \cite{KST} and $63/64$ of
\cite{Gebauer}. In fact, several of the open problems of the survey
of Gebauer, Moser, Scheder, and Welzl \cite{GMSW}
are centered around the understanding of this question.

In our main theorem we settle these questions from \cite{GMSW} and
determine the asymptotics of $f(k)$.
We show that the lower bound (\ref{kst}) can be strengthened by a
factor of $2$ and that this bound is tight.

\begin{thm} \label{main}
$$\left\lfloor\frac{2^{k+1}}{e(k+1)}\right\rfloor\le f(k)=\left(\frac{2}{e}+O\left(\frac{1}{\sqrt{k}}\right)\right)
\frac{2^{k}}{k}.$$
\end{thm}

For the upper bound, which constitutes the main contribution of our paper,
 we use the fundamental binary tree approach of \cite{Gebauer}.
The main novelty of our approach here is the development of a
suitable continuous setting for the construction of the
appropriate binary trees, which allows us to study the problem
via a differential equation. The solution of this differential equation
corresponds to our construction of the binary trees, which then
can be given completely discretely.

The lower bound is achieved via the lopsided version of the Lov\'asz Local
Lemma.
For the proof we set the values of the variables randomly and independently,
but not according to the uniform distribution. This seems reasonable to do as
the number of appearances of a variable $x_i$ in a CNF formula $F$ as a non-negated literal
could be significantly different from the number of clauses where $x_i$ appears
negated. It is even possible that a variable $x_i$ appears
negated in only a few, maybe even in just a single clause, in which case one tends
to think that it is reasonable to set this variable to {\tt true} with much larger
probability than setting it to {\tt false}. In fact, it is exactly the opposite we
will do. The {\em more} a variable appears in the clauses of $F$ as
non-negated, {\em the less likely} we will set it to {\tt true}.
The lower bound could also be derived from a theorem of
Berman, Karpinski and Scott \cite{BKS} tailored
to give good lower bounds on $f(k)$ for small values of $k$.
However the proof of \cite{BKS} contains a couple of inaccuracies obscuring
exactly this counterintuitive choice of the probabilities (in fact in \cite{BKS} the
probabilities are defined the opposite, the intuitive way).
Furthermore in \cite{BKS} the asymptotic behavior of the bound is not calculated,
since it was not believed to be optimal.
In Section \ref{sec:ProofOfLowerBound} we reproduce a simple
argument giving the asymptotics.

Since the (Lopsided) Lov\'asz Local Lemma was fully algorithmized by Moser and
Tardos \cite{MT} we now have that not only every $(k, s)$-CNF formula for
$s=\left\lfloor\frac{2^{k+1}}{e(k+1)}\right\rfloor$ has a satisfying
assignment but there is also an algorithm that {\em finds} such an assignment
in probabilistic polynomial time.
Moreover, for just a little bit larger value of the parameter $s$ one is not
likely to be able to
find a satisfying assignment efficiently, simply because already the decision
problem is NP-hard.

Our construction also shows
that the lower bound (\ref{lll}) on $l(k)$ is asymptotically tight.
\begin{thm}\label{thl}
$$l(k)=\left(\frac{1}{e}+O\left(\frac{1}{\sqrt{k}}\right)\right) 2^{k}.$$
\end{thm}

Theorem \ref{main} and Theorem \ref{thl} are another instances which show the tightness
of the  Lov\'asz Local Lemma. The first such example was given by Shearer \cite{Shearer}.

\subsection{$(k,d)$-trees}
The substantial part of the proofs of Theorems~\ref{main} and \ref{thl}---the upper bounds---as well as all our further results  depend on the construction of certain binary trees.

Throughout the paper, whenever mentioning {\em binary tree}s we always mean
{\em proper binary trees}, that is, rooted trees where every non-leaf node has exactly two children.
We say that a leaf $w$ of a tree $T$ is \emph{$k$-close} from a
vertex $v\in V(T)$ if $v$ is an ancestor of $w$ and $w$ has distance at most
$k$ from $v$. When $k$ is clear from the context we say {\em $w$ is visible
  from $v$} or {\em $v$ sees $w$}.

The concept of $(k,d)$-trees, introduced by Gebauer~\cite{Gebauer}, will be
our main tool in this paper. We call a (proper) binary tree $T$ a
\emph{$(k,d)$-tree}\footnote{To simplify our statements
we shifted the parameter $k$ in the definition compared to \cite{Gebauer}:
what is called a $(k,d)$-tree there is a $(k-1,d)$-tree in our new terminology.} if
\begin{itemize}
\item[(i)] every leaf has depth at least $k$ and
\item[(ii)] for every node $u$ of $T$ the number of $k$-close leaves from $u$ is at
  most $d$.
\end{itemize}
For a fixed $k$, we are interested in how low one can make $d$ in
a $(k,d)$-tree.
Essentially all our main results will be consequences of the construction of
$(k,d)$-trees with relatively small $d$.
We introduce $\f(k)$ to stand for the smallest
integer such that a $(k,\f(k))$-tree exists and determine $\f(k)$
asymptotically.\footnote{Note that in \cite{GMSW} the function
$\f$ is defined to be one less
than our definition (the maximum
$d$ such that {\em no} $(k,d)$-tree exist). This might be more appropriate if
$\f$ is only considered for its implications to SAT, but
in view of the many other applications we feel that our definition is more
natural.}

\begin{thm}\label{kd-trees}
$$\left\lfloor\frac{2^{k+1}}{e(k+1)}\right\rfloor<\f(k)=\left(\frac{2}{e}+O\left(\frac{1}{\sqrt{k}}\right)\right)
\frac{2^{k}}{k}.$$
\end{thm}

The construction of the trees providing the upper bound constitutes a large
portion of our paper. We devote quite a bit of effort (the entire Section~\ref{sec:Informal}) to describe the key informal ideas of the proof. That is, we formulate the construction process in a continuous setting and study its
progress with the help of a continuous two-variable function $F(t,x)$
defined by a certain differential equation.
It can be shown that
the integral $\int F(t,x)dx$ being large for some $t$ corresponds
to the construction process terminating with the desired $(k,d)$-tree.
Even though our treatment in Section~\ref{sec:Informal} will be highly informal (with simplifying assumptions and approximations) and is eventually not necessary for the formal 
proof, we find it important for a couple of reasons. On the one hand,
it provides the true motivation behind our formal discrete construction and
illustrates convincingly {\em why} this  construction (treated formally
in Section~\ref{sec:FormalProof}) {\em should} work.
Furthermore the continuous function $F(t,x)$, defined via the differential equation, is also helpful in studying the {\em size} of the constructed formulas.
This connection will be indicated in Section~\ref{sec:Outlook}.

\subsection{Formulas and the class MU(1)}

The function $f(k)$ is not known to be computable. In order to still be able
to upper bound its value, one tries to restrict to a smaller/simpler class of
formulas. When looking for unsatisfiable $(k,s)$-CNF formulas it is naturally
enough to consider {\em minimal unsatisfiable formulas}, i.e., unsatisfiable
CNF formulas that become satisfiable if we delete any one of their
clauses. The set of minimal unsatisfiable CNF formulas is denoted by MU.
As observed by Tarsi (cf.\ \cite{AharoniLinial}), all formulas in MU have more
clauses than variables, but some have only one more. The class of these
MU formulas, having one more clauses than variables, is denoted by MU(1).
This class has been widely studied
(see, e.g., \cite{AharoniLinial}, \cite{DDK}, \cite{KZ}, \cite{K}, \cite{S}).
Hoory and Szeider \cite{HSCompMU} considered the function $f_{1}(k)$, denoting
the largest integer such that no $(k, f_{1}(k))$-CNF formula is in MU(1), and
showed that $f_{1}(k)$ is computable.
Their computer search determined the values of $f_1(k)$ for
small $k$:
$f_1(5) = 7$, $f_1(6) = 11$, $f_1(7) = 17$, $f_1(8) = 29$, $f_1(9) = 51$.
Via the trivial inequality $f(k) \leq f_{1}(k)$, these are the best
known upper bounds on $f(k)$ in this range.
In contrast, even the value of $f(5)$ is not known. (For $k\le4$ it is
a doable exercise to
see that $f(k)=f_1(k)=k$.)

In \cite{Gebauer} $(k,d)$-trees were introduced to construct
unsatisfiable CNF formulas and
upper bound the functions $f$ and $l$. Since these formulas also
reside in the class MU(1), we can also use them to upper bound $f_1$. Most of
the content of the following statement appears already in Lemma~1.6 of
\cite{Gebauer}.

\begin{thm}[\cite{Gebauer}]\label{connection}
\begin{itemize}
\item[\rm(a)]$f(k)\le f_1(k)<\f(k)$
\item[\rm(b)] $l(k)<k\f(k-1)$
\end{itemize}
\end{thm}

For completeness and since part (a) is stated in a slightly weaker form in \cite{Gebauer},
we include the short proof in Section~\ref{subsection:formulas}.

It is an interesting open problem whether $f(k)=f_1(k)$ for every $k$.
Theorems~\ref{main}, \ref{kd-trees}, and \ref{connection} imply that $f(k)$ and
$f_1(k)$ are equal asymptotically: $f(k)=(1+o(1))f_1(k)$. More precisely, we have the following.

\begin{cor}
 $f_{1}(k) = \left(\frac{2}{e}+O\left(\frac{1}{\sqrt{k}}\right)\right)
\frac{2^{k}}{k}.$
\end{cor}
\noindent
Scheder \cite{Sched} showed that for almost disjoint $k$-CNF formulas
(i.e., CNF-formulas where any two clauses have at most one variable in common)
the two functions are not the same. That is, if $\tilde{f}(k)$ denotes the
maximum $s$ such that every almost disjoint $(k,s)$-CNF formula is
satisfiable, for $k$ large enough every unsatisfiable almost disjoint $(k,
\tilde{f}(k) + 1)$-CNF formula is outside of MU(1).

\subsection{The Neighborhood Conjecture}
A \emph{hypergraph} is a pair $(V,{\cal F})$, where $V$ is a finite set whose elements are
called \emph{vertices} and ${\cal F}$ is a family of subsets of $V$, called \emph{hyperedges}.
If it does not cause any confusion, we often refer to just ${\cal F}$ as the hypergraph.
A hypergraph is \emph{$k$-uniform} if every hyperedge contains exactly $k$ vertices.

A hypergraph is $2$-colorable if there is a coloring of the vertices with red and blue such that
no edge is {\em monochromatic}.

A standard application of the first moment method says that for any $k$-uniform hypergraph ${\cal F}$, we have
$$|{\cal F}| < 2^{k-1} \Rightarrow \mbox{${\cal F}$ is $2$-colorable}.$$
An important generalization of this implication was given by Erd\H os and Selfridge \cite{ESe}, which also
initiated the derandomization method of  conditional expectations.
Erd\H os and Selfridge formulated their result in the context of positional games.
Given a $k$-uniform hypergraph $\mathcal{F}$ on vertex set $V$,
players Maker and Breaker take turns in claiming
one previously unclaimed element of $V$, with Maker going first.
Maker wins if he claims all vertices of some hyperedge of $\mathcal{F}$, otherwise Breaker wins.

Since this is a finite perfect information game and the players have
complementary goals, either Maker has a winning strategy (that is, the
description of  the vertex to be claimed next by Maker in any imaginable game
scenario, such that at the end  he wins, no matter how Breaker plays)
or Breaker has a winning strategy. Which of this is the case depends solely on
${\cal F}$, hence it makes sense to call the hypergraph ${\cal F}$ {\em
Maker's win} or {\em Breaker's win}, respectively. Moreover, these games are monotone in the sense that adding hyperedges can only make it easier for Maker to win. In other words, if ${\cal F}$ is a Maker's win then ${\cal F}$ with additional hyperedges is naturally also a Maker's win.

The crucial connection between Maker/Breaker games to $2$-colorability is the
following:
\begin{equation}\label{eq:breaker-2-color}
\mbox{${\cal F}$ is Breaker's win $\Rightarrow$ ${\cal F}$ is $2$-colorable}.
\end{equation}
Indeed, if both players use Breaker's winning strategy,\footnote{The second
  player can use this strategy directly, but for the
  player going first one has to slightly modify it: start with an
  arbitrary move and then use the strategy, always maintaining to own one
  extra element of the board. It is easy to see
  that owning an extra element cannot be disadvantageous, so there exists a
  Breaker-win strategy for the first player, too.} by the end of the game they both win
as Breakers and hence create a $2$-colored vertex set $V$, where both colors
are represented in each hyperedge --- a proper $2$-coloring of ${\cal F}$.
Hence, the following theorem of Erd\H os and Selfridge \cite{ESe} is a generalization of the first moment result above.
$$|{\cal F}| < 2^{k-1} \Rightarrow \mbox{${\cal F}$ is a Breaker's win}.$$

As the Erd\H os-Selfridge Theorem can be considered the game-theoretic first moment method, the Neighborhood
Conjecture of J\'ozsef Beck (to be stated below) would be the game-theoretic Local Lemma.
Unlike the first moment method, the Local Lemma guarantees the $2$-colorability of hypergraphs based on some local condition
like the maximum degree of a vertex or an edge of the hypergraph.
The \emph{degree} $d(v)$ of a vertex $v\in V({\cal F})$ is the number of hyperedges of ${\cal F}$ containing $v$ and the
\emph{maximum degree} $\Delta(\mathcal{F})$ of $\mathcal{F}$ is the maximum degree of its vertices.
The \emph{neighborhood} $N(e)$ of a hyperedge $e$ is the set of hyperedges of $\mathcal{F}$ which intersect $e$, excluding $e$ itself, and the \emph{maximum neighborhood size} $\Delta ( L({\cal F}))$ of $\mathcal{F}$ is the maximum of $|N(e)|$
where $e$ runs over all hyperedges of $\mathcal{F}$. ($L({\cal F})$ denotes the {\em line-graph} of ${\cal F}$).
Simple applications of the Local Lemma show that,
\begin{align}
 \Delta( L({\cal F})) \leq \frac{2^{k-1}}{e} -1 & \Rightarrow \mbox{${\cal F}$ is $2$-colorable,} \label{lll1}\\
\Delta ( {\cal F}) \leq \frac{2^{k-1}}{ek}  & \Rightarrow \mbox{${\cal F}$ is $2$-colorable}. \label{lll2}
\end{align}

The Neighborhood Conjecture in its strongest form \cite[Open Problem
9.1]{Beckbook}
was suggesting the far-reaching generalization
that already when $\Delta ( L({\cal F})) < 2^{k-1} -1$,
${\cal F}$ should be a Breaker's win.  This was motivated by the construction of
Erd\H{o}s and Selfridge of a $k$-uniform Maker's win hypergraph $(V,\mathcal{G})$ with
$|{\cal G}|= 2 ^{k - 1}$, showing the tightness of their theorem.
The maximum neighborhood size of this hypergraph is $2^{k - 1}-1$ (every pair of edges intersects) and no better construction was known
until Gebauer~\cite{Gebauer} disproved the conjecture using her $(k,d)$-tree approach.
She constructed Maker's win hypergraphs ${\cal F}$ and ${\cal H}$ with $\Delta (L({\cal F})) = 0.75\cdot 2^{k-1}$ and
$\Delta ({\cal H}) \leq \frac{63}{128}\frac{2^{k}}{k}$, respectively.
Our $(k,d)$-trees from Theorem \ref{main} will imply the following somewhat improved bounds.
\begin{thm}\label{maker-breaker} For every integer $k \geq 3$ there exists
a Maker's win $k$-uniform hypergraph ${\cal H}$ such that
\begin{itemize}
\item[$(i)$] $\Delta ( {\cal H}) \le
  2\f(k-2)=\left(1+O\left(\frac{1}{\sqrt{k}}\right)\right) \frac{2^{k}}{ek}$
  \, and
\item[$(ii)$] $\Delta ( L ({\cal H})) =(k-1)\f(k-2)=
  \left(1+O\left(\frac{1}{\sqrt{k}}\right)\right) \frac{2^{k-1}}{e}$.
\end{itemize}
\end{thm}
Note that the bound in part $(ii)$ asymptotically coincides with the one
given by the Local Lemma in \eqref{lll1} for 2-colorability,
while part $(i)$ is still a factor $2$ away from \eqref{lll2}.
Note furthermore that the bounds in \eqref{lll1} and \eqref{lll2} are not
optimal. More elaborate methods of Radhakrishnan and Srinivasan \cite{RS}
show that for a small enough constant $c>0$ any $k$-uniform hypergraph ${\cal
  F}$ with $\Delta (L({\cal F})) \le c2^k \sqrt{k/\log k})$ is, in fact,
2-colorable. Part $(ii)$ of Theorem~\ref{maker-breaker} establishes that no 
game-theoretic $2$-colorability condition can exist 
beyond the asymptotic Local Lemma bound.

In Section~\ref{newest} we return to the game theoretic questions
discussed here and elaborate on what are the winning strategies of the
players that we can analyze (the so called {\em pairing strategies}),
and why they are insufficient to resolve weaker forms of the Neighborhood Conjecture
in either direction.

\subsection{European Tenure Game}

The (usual) Tenure Game (introduced by Joel Spencer~\cite{Spencer})
is a perfect information game between two players: the (good) chairman of the department,
and the (vicious) dean of the school. The department has $d$
non-tenured faculty and the goal of the chairman is to promote
(at least) one of them to tenure, the dean tries to prevent this. Each
non-tenured faculty is at one of $k$ pre-tenured rungs, denoted
by the integers $1,\ldots ,k$. A non-tenured faculty becomes tenured if she has rung $k$ and is promoted.
The procedure of the game is the following. Once each year, the chairman proposes to the dean
a subset $S$ of the non-tenured faculty to be promoted by one rung. The dean has two choices: either he accepts the suggestion
of the chairman, promotes everybody in $S$ by one rung and fires everybody
else, or he does the complete opposite of
the chairman's proposal (also typical dean-behavior): fires everybody in $S$ and promotes everybody else by one rung.
This game obviously ends after at most $k$ years. The game analysis is very
simple, see \cite{Spencer}.

In the European Tenure Game (introduced by Benjamin Doerr~\cite{Doerr04})
the rules are modified, so the non-promoted part of the non-tenured faculty is
not fired, rather demoted back to rung 1. An equivalent, but perhaps more
realistic scenario is that the non-promoted faculty is fired but the
department hires new people at the lowest rung to fill the tenure-track
positions vacated by those fired.
For simplicity we assume that all non-tenured faculty are at the lowest rung
in the beginning of the game and would like to know what combinations of $k$
and $d$ allow for the chairman to eventually give tenure for somebody when
playing against any (vicious and clever) dean. For fixed $d$ let $v_d$ stand for
the largest number $k$ of rungs such that this is possible.

Doerr~\cite{Doerr04} showed that
$$\lfloor \log d +\log\log d +o(1)\rfloor \leq v_d \leq \lfloor \log d +\log\log d + 1.73 +o(1)\rfloor.$$
It turns out that the game is equivalent to $(k,d)$-trees, hence using
Theorem~\ref{kd-trees} we can give a precise answer, even in
the additive constant, which turns out be $\log_2 e -1 \approx 0.442695$.
\begin{thm} \label{european}
The chairman wins the European Tenure Game with $d$ faculty and $k$ rungs if
and only if there exists a $(k,d)$-tree.  In particular,
$v_d=\max\{k\mid\f(k)\le d\}$ and we have
$$v_d = \lfloor \log d +\log\log d + \log e  - 1 + o(1)\rfloor.$$
\end{thm}

\subsection{Searching with lies}

In a liar game the first player, called Chooser, thinks of a member $x$ of an agreed upon $N$ element
set $H$ and the second player, called Guesser, tries to figure it out by Yes/No questions of the sort
``Is $x\in S$?'', where $S$ is a subset of $H$ picked by Guesser.
This is not difficult if Chooser is always required
to tell the truth, but usually  Chooser is allowed to lie. However for Guesser to have a chance to be successful, the lies of Chooser
have to come in some controlled fashion. The most prominent of these
restrictions allows Chooser to lie at most $k$ times and asks for the smallest
number $q(N,k)$ of questions that allows
Guesser to figure out the answer.
This is also called Ulam's problem for binary
search with $k$ lies.
For an exhaustive description of various other lie-controls, see the survey of Pelc~\cite{Pelc}.

One of the problems in the 2012 International Mathematics Olympiad
was a variant of the liar game. Instead of limiting the total number of lies,
in the IMO problem the number of consecutive lies was limited. This
fits into the framework of Section 5.1.3 in Pelc's survey~\cite{Pelc}.
This restriction on the lies is not enough for Guesser to find the value $x$ with
certainty, but he is able to narrow down the set of possibilities. The IMO
problem asks for certain estimates on how small Guesser can eventually make this
set.
This problem was also the topic of the Minipolymath4 project research thread
\cite{polymath}.

It turns out that this question can also be expressed in terms of existence of
$(k,d)$-trees.

\begin{thm}\label{IMO}
Let $N>d$ and $k$ be positive integers. Assume Chooser and Guesser play the guessing game
in which Chooser thinks of an element $x$ of an agreed upon set $H$ of size $N$ and
then answers an arbitrary number of Guesser's questions of the form ``Is $x\in
S$?''. Assume further that Chooser is allowed to lie, but never to $k$
consecutive questions. Then Guesser can guarantee to narrow the number of possibilities
for $x$ with his questions to at most $d$ distinct values if and only if a
$(k,d+1)$-tree exists, that is, if $d<\f(k)$.
\end{thm}

\subsection{Notation} Throughout this paper, $\log$ denotes the binary
logarithm. We use $\N$ to denote the set of natural numbers including $0$.

As we have mentioned, by a \emph{binary tree} we always mean a rooted tree
where every node has either two or no children. The distance between two vertices $u$, $v$ in a binary tree is the number of edges in the unique path from $u$ to $v$.
The {\em depth} of a vertex is its distance from the root, and the {\sl height} of a binary tree is the maximum depth of a vertex.

\subsection{Organization of this paper}
In Section~\ref{sec:applications} we derive all applications of the upper
bound in Theorem~\ref{kd-trees}.
This includes the upper bounds of Theorems \ref{main} and \ref{thl} through
proving Theorem~\ref{connection}, as well as
the proofs of
Theorems~\ref{maker-breaker}, \ref{european} and
\ref{IMO}.
In Section~\ref{sec:FormalDef} we give the basic definitions and simple propositions required for the proof of Theorem~\ref{kd-trees}.
In Section \ref{sec:Informal} we sketch the main informal ideas behind our tree-construction
including a rough description of our approach, and how it translates to
solving a differential equation. This is in the background of
our actual formal constructions for the proof of Theorems~\ref{kd-trees}, which then can be given completely discretely.
The formal construction is the subject of Section \ref{sec:FormalProof}. The lower bound of
Theorem \ref{main}
 (implying the lower bound for Theorem~\ref{kd-trees})
is shown in
Section \ref{sec:ProofOfLowerBound}.
In Section \ref{sec:Outlook} we give an outlook and pose some open problems.


\section{Applying $(k,d)$-trees} \label{sec:applications}
In this section we apply $(k,d)$-trees and the upper bound in Theorem~\ref{kd-trees} to prove the upper bounds of
Theorems \ref{main} and \ref{thl}, as well as
Theorems~\ref{connection}, \ref{maker-breaker}, \ref{european}, and \ref{IMO}.
Some of these connections, sometimes in disguise, were already pointed out in \cite{Gebauer}.

\subsection{Formulas}\label{subsection:formulas}

In this subsection we give a proof of Theorem~\ref{connection}, which, together with the upper bound in Theorem~\ref{kd-trees},
readily implies the upper bounds in Theorems~\ref{main} and \ref{thl}.

For every binary tree $T$ (recall that we only consider proper binary trees) which has all of its leaves at depth at least $k$,
one can construct a $k$-CNF formula $F_k(T)$ as follows.
For every non-leaf node $v\in V(T)$ we create a variable
$x_{v}$ and label one of its children with the literal $x_v$ and the other with
$\bar{x}_{v}$.
We do not label the root.
With every leaf $w\in V(T)$
we associate a clause $C_w$, which is the conjunction of the first $k$ labels encountered when
walking along the path from $w$ towards the root (including the one at $w$).
The disjunction of the clauses $C_{w}$ for all leaves $w$ of $T$ constitutes
the formula $F_k(T)$.
\begin{obse} \label{obse:Funsat}
$F_k(T)$ is unsatisfiable.
\end{obse}
\begin{proof}
Any assignment $\alpha$ of the variables defines a path from
the root to some leaf $w$ by always proceeding to the unique child
whose label is mapped to {\tt false} by $\alpha$.
Then $C_w$ is violated by $\alpha$.
\end{proof}

\begin{proof}[Proof of Theorem~\ref{connection}]
Consider a $(k,\f(k))$-tree $T$ and the corresponding $k$-CNF
formula $F=F_k(T)$. $F$ is unsatisfiable by Observation~\ref{obse:Funsat}. The
variables of this formula are the vertex-labels of
$T$. The variable $x_v$ corresponding to a vertex $v$ appears in the
clause $C_w$ if and only if $w$ is
$k$-close from $v$. Thus each variable appears at most $\f(k)$ times. This
makes $F$ an unsatisfiable $(k,\f(k))$-CNF, proving $f(k)<\f(k)$. If
$F$ is in MU(1), then we further have $f_1(k)<\f(k)$. The number of
clauses in $F$ is the number of leaves in $T$, the number of variables in $F$
is the number of non-leaf vertices in $T$, so we have exactly one more clauses
than variables. But $F$ is not a {\em minimal} unsatisfiable formula in
general. Fortunately it is one, if each variable appears in $F$ both in negated and
in non-negated forms, see \cite{DDK}. This will be the
case if we pick $T$ to be a $(k,\f(k))$-tree
which is minimal with respect to containment.
Indeed, if a literal
associated to a vertex $v$ does not appear in any of the clauses, then the
subtree of $T$ rooted at $v$ is a $(k,\f(k))$-tree.
This finishes the proof of part (a) of the theorem.

Clearly, the neighborhood of any clause in a $(k,d)$-CNF formula is of size
at most
$k(d-1)$, but this bound is too rough to prove part (b) with. We start by
picking a $(k-1,\f(k-1))$-tree $T'$ with each leaf having depth at least
$k$. Such a tree can be constructed by taking two copies of an arbitrary
$(k-1,\f(k-1))$-tree and connecting their roots to a new root vertex. Now
$F'=F_k(T')$ is an unsatisfiable $k$-CNF by
Observation~\ref{obse:Funsat}. The advantage of this construction is that now
each {\em literal} appears at most $\f(k-1)$ times (as opposed to only having
a bound on the multiplicity of each {\em variable} in $T$). Note that if two
distinct clauses in $F'$ share a common variable, then there is a variable
that appears negated in one of them and non-negated in the other. This implies
that every clause intersects at most $k\f(k-1)$ other clauses as needed.
\end{proof}

Note that $F'$ in the proof above can also be chosen to be in MU(1). This is
the case if $T'$ is obtained from a minimal $(k-1,\f(k-1))$-tree by doubling
it.

Now Theorem~\ref{kd-trees} together with part (a) of
Theorem~\ref{connection} implies the upper bound in Theorem~\ref{main},
while together with part (b) of Theorem~\ref{connection}
it implies the upper bound in Theorem~\ref{thl}. We also have an implication in
the reverse direction: the lower bound of Theorem~\ref{main} implies the lower
bound of Theorem~\ref{kd-trees} using Theorem~\ref{connection}(a).

\subsection{The Neighborhood Conjecture} \label{sec:improvingneighborhoodconjecture}

In this subsection we prove Theorem~\ref{maker-breaker}.

We start with a few simple observations. One can associate a hypergraph ${\mathcal H}(F)$
to any CNF formula $F$ by taking all the literals in $F$ as vertices, and
considering the clauses of $F$ (or rather the set of literals which the
clause is the disjunction of) as the hyperedges.
\begin{obse}\label{hf}
If $F$ is unsatisfiable, then ${\mathcal H}(F)$ is Maker's win.
\end{obse}
\begin{proof}
As we have seen before, going second cannot help Maker, so we assume that Breaker starts. Now whenever Breaker picks a literal $u$, Maker immediately picks $\bar{u}$ in the following move.
At the end, consider
the evaluation of the variables setting all the literals that Breaker has
to {\tt true}.
As $F$ is unsatisfiable this evaluation falsifies $F$ and therefore it
violates one of the clauses, giving Maker his win at the hyperedge
corresponding to this clause.
\end{proof}

As we have seen in the proof of Theorem~\ref{connection} there exist a
$(k-1,\f(k-1))$-tree $T'$ with all leaves at depth at least $k$ and $F_k(T')$ is
an unsatisfiable $k$-CNF where all literals appear in at most $\f(k-1)$
clauses. By Observation~\ref{hf} this makes ${\mathcal H}(F_k(T'))$ a Maker's
win $k$-uniform
hypergraph with maximum degree at most $\f(k-1)$. It is easy to see that
$\f(k-1)\le2\f(k-2)$, so this bound is better than the one claimed in part (i)
of Corollary~\ref{maker-breaker}. But in order to get part (ii) as well we have to
take a slightly different approach.

\begin{proof}[Proof of Theorem~\ref{maker-breaker}]
Let $T$ be a $(k-2,\f(k-2))$-tree, with all leaves in depth at least $k$. One
can build such a tree from four copies of an arbitrary
$(k-2,\f(k-2))$-tree. In each clause of $F=F_k(T)$ the literals are associated
to vertices in $T$. We distinguish the {\em leading} literal in each clause to
be the one associated to the vertex closest to the root. A clause $C_w$
contains the literal associated to a vertex $v$ in {\em non-leading} position
exactly when the leaf $w$ is $(k-2)$-close from $v$. Thus any literal
appears in at most $\f(k-2)$ clauses in non-leading position. A literal
associated to a leaf of $T$ appears in a single clause of $F$. Any other
literal $\ell$ appears in at most $2\f(k-2)$ clauses as any clause
containing $\ell$ must also contain
the literal associated to one of the two children of the vertex of $\ell$ and these are not in
leading position. This makes $\mathcal H={\mathcal H}(F)$ a $k$-uniform
hypergraph with
maximum degree $\Delta(\mathcal H)\le2\f(k-2)$. Observations~\ref{obse:Funsat}
and \ref{hf} ensure that $\mathcal H$ is a Maker's win hypergraph and
this concludes the proof of part (i).

As we have observed earlier each pair of distinct intersecting clauses of $F$
share a variable that appears negated in one of them and in non-negated form
in the other. Observe further that when two distinct clauses share a literal
they also share a variable that appears in opposite form as non-leading
literals in the two clauses. A clause $C$ has $k-1$ non-leading literals and
the opposite form of each is contained in non-leading position in at most $\f(k-2)$ clauses. This gives us the bound
stated in part (ii) on the number of clauses sharing a literal with $C$.
\end{proof}

\subsection{The European Tenure Game}
In this subsection we prove Theorem~\ref{european}, namely the connection to
$(k,d)$-trees. The formula estimating $v_d$ will then follow from
Theorem~\ref{kd-trees}. We start with a special labeling of $(k,d)$-trees.

\begin{prop}\label{labeling}
Let $T$ be a $(k,d)$-tree and let $L$ be its set of leaves. There exists a
labeling $i:L\to\{1,\ldots,d\}$ such that for every vertex $v\in V(T)$ all
the $k$-close leaves from $v$ have distinct labels.
\end{prop}
\begin{proof}
We define the labels of the leaves one by one.
We process the vertices of $T$  according to a Breadth First Search, starting
at the root. When processing a vertex $v$ we label the still unlabeled leaves
that are $k$-close from $v$ making sure they all receive different
labels. This is possible, because the ones already labeled are all
$k$-close from the parent of $v$, so must have received different
labels. We have enough labels left because the total number of leaves visible
from $v$ is at most $d$. After processing all vertices of $T$ our labeling is
complete and satisfies the requirement.
\end{proof}

\begin{proof}[Proof of Theorem~\ref{european}]
Suppose first that there is a $(k,d)$-tree $T$ and let us give a winning strategy to the
chairman. We start by labeling the leaves of $T$ with the $d$
non-tenured faculty,  according to
Proposition~\ref{labeling}.
The chairman is placed on the root of $T$ and during the game he will move
along a path from the root to one of the leaves, in each round
proceeding to one of the children of its current position.
To which leaf he arrives at depends on the answers of the
dean. When standing at a non-leaf vertex $v$
with left-child $v_1$ and right-child $v_2$, the chairman proposes to promote the subset $S$ of the
faculty containing the labels of the leaves that are $(k-1)$-close from
$v_1$. If the dean accepts his proposal he moves to $v_1$, otherwise he moves
to $v_2$. The game stops when a leaf is reached. We claim that the label $P$ of
this leaf is promoted to tenure, hence the chairman has won.

Note that by part (i) of the definition of a $(k,d)$-tree the game lasts for
at least $k$ rounds.
We will show that in each of the last $k$ rounds $P$ was promoted.
Indeed, if the chairman moved to the left in one of these rounds then $P$ was
proposed for promotion and the dean accepted it. However, if the chairman moved
to the right, then $P$ could not be proposed for promotion by
the condition of the labeling of Proposition~\ref{labeling}, but the dean
reversed the proposal and hence
$P$ was promoted in these rounds as well.

For the other direction, we are given a winning strategy for the
chairman. This strategy specifies at any point of the game which subset of the
faculty the chairman should propose for promotion unless a member of the
faculty is already tenured, at which time the game stops. In building the game tree
we disregard the subsets but pay close attention to when the game stops. In
particular, each vertex corresponds to a position of the game with the root
corresponding to the initial position. If a vertex $v$ corresponds to a
position where the game stops we make $v$ a leaf and label it with one of the
faculty members that has just been tenured. Otherwise $v$ has two children,
one corresponding to each of the two possible answers of the dean.

Clearly, this is a (proper) binary tree.
We claim that it is a $(k,d)$-tree. Note that
in order for somebody get tenured she has to be promoted in $k$ consecutive
rounds, so all leaves are at depth at least $k$, as required.

To prove that no vertex sees more than $d$ leaves we prove that all leaves
$k$-close from the same vertex have distinct labels. Indeed, if a leaf $w$
is $k$-close from a vertex $v$ and $w$ is labeled by faculty Frank, then
Frank had to be promoted in all rounds in the game from the position
corresponding to $v$ till the position corresponding to $w$. But Frank is
promoted in exactly one of the two cases depending on the dean's answer, so this
condition determines a unique path in our tree from $v$ making $w$ the unique
leaf labeled Frank that is $k$-close from $v$.
\end{proof}

\subsection{Searching with lies}

\begin{proof}[Proof of Theorem~\ref{IMO}]
The first observation is that the game with parameters $N>d$ and $k$ is won by
Guesser if and only if Guesser wins with parameters $N'=d+1$, $d$ and $k$.

One direction is trivial, if $N$ is decreased it cannot hurt Guesser's chances. For
the other direction assume $B$ has a winning strategy for the $N'=d+1$
case. This means that given $d+1$ possible values for $x$ he can eliminate one
of them. Now if he has more possibilities for the value of $x$ he can
concentrate on $d+1$ of them and ask questions till he eliminates one of these
possibilities. He can repeat this process, always reducing the number of
possible values to $x$ till this number goes below $d+1$ --- but by then he
has won.

We can therefore concentrate to the case $N=d+1$. We claim that this
is equivalent to the European tenure game with $k$ non-tenured rungs and $d+1$
non-tenured faculty, thus
Theorem~\ref{IMO} follows from Theorem~\ref{european}. Indeed, Chooser corresponds
to the dean, Guesser corresponds to the chairman, the $d+1$
possible values of $x$ correspond to the non-tenured faculty members. Guesser asking
whether $x\in S$ holds corresponds to the
chairman proposing the set $S$ for promotion, a ``no'' answer corresponds to
the dean accepting the proposal, while a ``yes'' answer corresponds to the
dean reversing it. At any given time the rung
of the faculty member corresponding to possible value $v$ in the liar
game is $i+1$, where $i$ is the largest value such that the last $i$ questions
of Guesser were answered by Chooser in a way that would be false if $x=v$. Thus, a win for the chairman
(tenuring a faculty member, i.e. promoting him $k$ consecutive times)
exactly corresponds to Guesser answering $k$ consecutive times in such a way
that would be false if $x=v$. This makes $x=v$ impossible according to
the rules of the liar game, so Guesser can eliminate $v$ and win.
\end{proof}


\section{Formal Definitions and Basic Statements}
\label{sec:FormalDef}

\subsection{Vectors and constructibility}
Given a node $v$ in a tree $T$, it is important to count the leaf-descendants
of $v$ in distance $i$ for $i\le d$. We say that a non-negative integer vector
$(x_{0}, x_{1}, \ldots, x_{k})$ is a {\em leaf-vector} for
$v$ if $v$ has at most $x_{i}$ leaf-descendants in distance $i$ for each $0\le
i\le k$.
E.g., the vector $(1, 0, \ldots, 0)$ is a leaf-vector for any leaf, while
for the root $v$ of a full binary tree of height $l \leq k$
we have $(0, 0, \ldots, 0, 2^{l}, 0, \ldots, 0)$ as its smallest leaf-vector.
We set $|\vec{x}| := \sum_{i=0}^{k} x_{i}$.
By definition, every node $v$ of a $(k,d)$-tree has
a leaf-vector $\vec x$ with $|\vec{x}|  \leq d$.

For some vector $\vec{x} \in \N^{k+1} $
we define a \emph{$(k,d,\vec{x})$-tree} to be a tree where
\begin{itemize}
\item[(i)] $\vec x$ is a leaf-vector for the root, and
\item[(ii)] each vertex has at most $d$ leaves that are $k$-close.
\end{itemize}
For example,
a tree consisting of a parent with two children is a
 $(k,d,(0,2,0\ldots,0))$-tree for any $k\geq 1$ and $d\geq 2$.

We say that a vector $\vec{x}\in\N^{k+1}$ is \emph{$(k,d)$-constructible}
(or \emph{constructible} if $k$ and $d$ are clear from the context),
if a $(k,d,\vec{x})$-tree exists.
E.g., $(1, 0, \ldots , 0)$, or more generally
$(\underbrace{0, 0, \ldots, 0}_{l}, 2^{l}, 0, \ldots, 0)$ are
$(k,d)$-constructible as long as $2^l\leq d$.

\begin{obse} \label{obse:constructibility}
There exists a $(k,d)$-tree if and only if
the vector $(0, \ldots, 0, d)$ is $(k,d)$-constructible.
\end{obse}
\begin{proof}
The vector $(0, \ldots, 0, d)$ is $(k,d)$-constructible if and only
if a $(k,d,(0,\ldots,0,d))$-tree exists.
It is easy to see that the definitions of $(k,d)$-tree and $(k,d,(0,\ldots,0,d))$-tree are equivalent.
Indeed, condition (ii) is literally the same for both, while condition (i)
states in both cases that there is no leaf $k-1$-close to the root. The only
difference is that condition (i) for a $(k,d,(0,\dots,0,d))$-tree also states
that there are at most $d$ leaves $k$-close to the root, but this also follows
from (ii).
\end{proof}

The next observation will be our main tool to capture how leaf-vectors change
as we pass
from a parent to its children.
\begin{obse}\label{obse:splitting}
If $\vec x'=(x'_0,x'_1,\ldots , x'_k)$ and $\vec x''=(x''_0,x''_1,\ldots , x''_k)$ are
$(k,d)$-constructible and $|\vec x|\le d$ for $\vec x=(0,x'_0+x''_0, x'_1+x''_1, \ldots
, x'_{k-1}+x''_{k-1})$, then $x$ is also $(k,d)$-constructible.
\end{obse}
\begin{proof}
Let $T'$ be a $(k,d, \vec x')$-tree with root $r'$ and $T''$ a $(k,d,\vec x'')$-tree
with root $r''$.
We create the tree $T$ by adding a new root vertex $r$ to the disjoint union of $T'$ and $T''$ and attaching it to both $r'$ and $r''$.
This tree is a $(k,d,\vec x)$-tree. Indeed, the leaf-descendants of $r$ at distance
$i$ are exactly leaf-descendants of either $r'$ or $r''$ at distance $i-1$,
hence $x$ is a leaf-vector for $r$. We also have to check that no
vertex has more than $d$ leaves $k$-close. This holds for the
vertices of $T'$ and $T''$ and is ensured by our assumption $|\vec x|\le d$ for the
root $r$.
\end{proof}

For a vector $\vec{x} = (x_{0}, \ldots, x_{k})$  we define its \emph{weight} $w(\vec{x})$
to be $\sum_{i = 0}^{k} x_{i}/2^{i}$.
The next lemma
gives a useful sufficient condition for
the constructibility of a vector.
\begin{lemma} \label{lemma:payofflemma}
 Let $\vec{x}\in\N^{k+1}$ with $|\vec x|\le d$.
If $w(\vec{x}) \geq 1$ then $\vec{x}$ is $(k,d)$-constructible.
\end{lemma}
\noindent
We note that Lemma \ref{lemma:payofflemma} is a reformulation of Kraft's
inequality.
For completeness we give a direct proof here.

\begin{proof}[Proof of Lemma  \ref{lemma:payofflemma}]
We build a binary tree starting with the root and adding the levels one
by one. As long as $\sum_{j = 0}^{i} \frac{x_{j}}{2^{j}} < 1$, we select a set
of $x_{i}$ vertices from the vertices on level $i$
and let them be leaves. We construct the $(i+1)$th level by adding two
children to the remaining
$2^{i}(1 - \sum_{j = 0}^{i} \frac{x_{j}}{2^{j}})$ vertices on level $i$.
At the first level $\ell \leq k$ where $\sum_{j = 0}^{\ell} \frac{x_{j}}{2^{j}} \geq 1$ we
mark all vertices as leaves and stop the construction
of the tree. The total number of leaves is at most $\sum_{j = 0}^{\ell}
x_{j} \leq |\vec x| \leq d$
and the number of leaves at distance $j$ from the root is at most $x_{j}$,
so the constructed tree is a $(k,d,\vec{x})$-tree.
\end{proof}

\begin{figure} [tbp]
\centering
\includegraphics[width=0.98\textwidth]{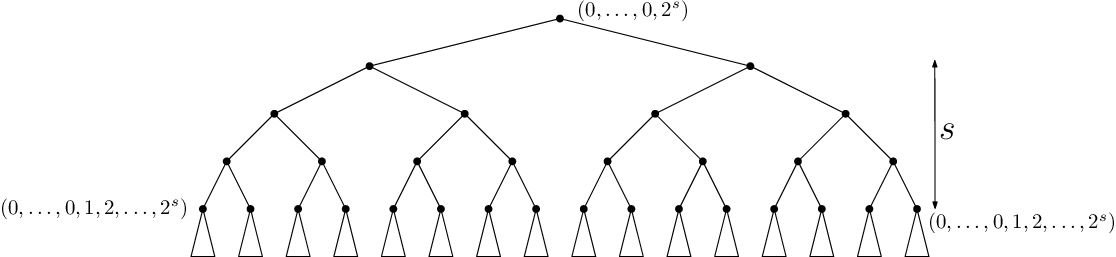}
\caption{Construction in \cite{Gebauer}: attaching a $(k,d, (0, \ldots, 0, 1,
  2, \ldots , 2^s))$-tree to every leaf of a full binary tree of height $s$
  gives a $(k,d)$-tree.} \label{fig:IlluWeakerBound}
\end{figure}

The main result of \cite{Gebauer} is the construction
of
$(k,d)$-trees with $d=\Theta (\frac{2^k}{k})$.
This argument is now streamlined via Lemma~\ref{lemma:payofflemma}.
Indeed, the $(k,d)$-constructibility
of the vector $\vec{v} = (0, \ldots, 0, 1, 2, 4,\ldots , 2^s)$ is an immediate
consequence of
Lemma~\ref{lemma:payofflemma}, provided $1\leq w(\vec{v}) = \sum_{i=k-s}^k
2^{i - k +s}/2^i=(s+1)2^{s-k} $
and $d\geq |\vec{v}| = \sum_{i=0}^k v_i = 2^{s+1} -1$.
Setting $s=k-\lfloor\log(k-\log k)\rfloor$ and $d=2^{s+1} - 1$ allows both
inequalities to hold.
Then by repeated application of Observation~\ref{obse:splitting} with $\vec x'=\vec x''$
we obtain
the constructibility of $(0,\ldots, 0,2,4\ldots , 2^s)$, $(0,\ldots ,0, 4,
\ldots , 2^s)$, etc.,
and finally the constructibility of $(0,\ldots , 0, 2^s)$.
This directly implies the existence of a $(k,d)$-tree by
Observation~\ref{obse:constructibility}. Note that $d=(2+o(1))\frac{2^{k}}k$
for infinitely many $k$ including
$k=2^t+t+1$ for any $t$.
Figure \ref{fig:IlluWeakerBound} shows an illustration.

\paragraph{Proof strategy.}

 In Section~\ref{sec:FormalProof} we will construct our $(k,d)$-tree starting
 with the root, from top to bottom.
When considering  some vertex $v$ it will be assigned a leaf-vector
$\vec{\ell}_v$. At this moment $v$ itself is a leaf in the partly
constructed tree, so one should consider $\vec{\ell}_v$ just as a promise:  for each
$i=0,1\ldots, k$ the vertex $v$ promises to have at most
$(\vec{\ell}_v)_i$ leaf-descendants at distance $i$
when the $(k,d)$-tree is fully constructed.

We start with the root with a leaf-vector $(0,\ldots , 0,d)$.
At each step we have to consider the vertices $v$ that are currently leaves but
promise not to be leaves: i.e., having a leaf-vector $\vec x$ with
$(\vec x)_0=0$. For such a vertex $v$ we add two children and associate
leaf-vectors $\vec x'$ and $\vec x''$ to them. According to
Observation~\ref{obse:splitting} we have to split the coordinates of $\vec x$
observing $x'_{i-1} +x''_{i-1} = x_i$ for $1\le i\le k$ and then we can decide
about the last coordinates $x'_k$ and $x''_k$ almost freely,
though we must respect the bounds $|\vec x'|\le d$ and $|\vec x''|\le d$.

We do not have to worry about nodes $v$ with a leaf-vector $\vec x$ satisfying
$w(\vec x) \geq 1$: Lemma~\ref{lemma:payofflemma}
ensures that $\vec x$ is constructible. Making $v$ the root of a $(k,d,\vec
x)$-tree we ensure that $v$ keeps its promise.

\medskip

The proof of Theorem~\ref{kd-trees}
contains several technically involved arguments.
In the next section we sketch
the main informal ideas behind the construction. Even though
the final construction can be formulated without mentioning the underlying
continuous context, we feel that an informal description greatly helps
in motivating it. Furthermore any future attempt to obtain more precise bounds 
on $f(k)$ most likely must consider the limiting continuous setup.
However, a reader in a hurry for a precise argument is
encouraged to skip right ahead to Section~\ref{sec:FormalProof}: 
the next section is not necessary for the formal understanding of that proof.


\section{Informal Continuous Construction}
\label{sec:Informal}

\subsection{Operations on leaf-vectors}
By the argument at the end of the last section all we care about from
now on
are leaf-vectors and how we split them up between the
two children, such that eventually all leaf-vectors have weight at least $1$.

We will consider two fundamentally different ways a parent
vertex $v$ with leaf-vector $(x_{0}, x_{1}, \ldots, x_{k})$
can split up its leaf-vector (in effect: its allotted number of $k$-close
leaves) between its children.
In the {\em fair split} both children get the same vector.
In this case the children can even get a last coordinate of at most $d/2$ and
their coordinate sum would still be at most $d$.
In the simplest case of the {\em piecewise split} the left child gets all the
leaves that are $t$-close
and the right child gets the $k$-close leaves whose distance is more than
$t$. In other words all the non-zero coordinates in the leaf-vector of the
left child will
be on the left of the non-zero coordinates of the leaf-vector of the right
child. For simplicity we keep the last coordinate of the leaf-vectors of the
children $0$. In the general case of the piecewise split we split a
leaf-vector to many vectors, one inheriting all $t$-close leaves, while the
others split the farther leaf-descendants evenly.

In the following informal description we will split leaf-vectors and divide
their coordinates freely, not caring about
divisibility. Dealing with rounding is one of the issues we leave to the
formal argument.

\paragraph{Fair Split.}
The leaf-vector $\vec x$ of the parent node $v$ is split evenly between its children $v'$ and
$v''$. Furthermore their last coordinate is $d/2$. That is,
$$\vec x'=\vec x''= E(\vec x) := (x_1/2,x_2/2,\ldots , x_k/2, d/2).$$
By $m$ repeated applications of the fair split we obtain  the leaf-vector
$$E^m(x) := \left( \frac{x_m}{2^m}, \frac{x_{m+1}}{2^m}, \ldots , \frac{x_k}{2^m}, \frac{d}{2^m}, \frac{d}{2^{m-1}}, \ldots , \frac{d}{2}\right).$$
After the $m$-times iterated fair split the leaf-vectors of all $2^m$ leaves
of the full binary tree so obtained are equal. After this operation it
is sufficient to ensure the constructibility of this single leaf-vector.
\medskip

In the previous section the iterated fair split was used on the leaf-vector
$(0,\ldots,2^s)$ to obtain the
leaf-vector $(0,\ldots , 0,1,2,4, \ldots , 2^s)$. The constructibility of the
latter vector was ensured by Lemma~\ref{lemma:payofflemma}.
The result obtained there is the best one can do using only fair splits and
Lemma~\ref{lemma:payofflemma} and is a factor $\frac{1}{e}$ away from
our goal. In order to improve we will also use the piecewise
splitting of the leaf-vectors, where the $l=1$ case of the piecewise split
can be thought of as a sort of complete opposite of the fair split.
\paragraph{Piecewise Split}
This split has two parameters, $r$ and $l$ with $1\le l\le r\le k$. Piecewise
split of a node $v$ with leaf-vector $\vec x=(x_0,\ldots,x_k)$ is
similar to the $l$-times iterated fair split in that we insert a depth $l$
full binary tree under $v$. But instead of assigning the
same leaf-vector to all $2^l$ leaves of this binary tree we assign a different
leaf-vector $\vec x'$ to one leaf and the same leaf-vector $\vec x''$ to
the remaining $2^l-1$ leaves of this full binary tree. We call the node with
leaf-vector $x'$ the left-descendant and the ones with leaf-vector $\vec x''$
the right-descendants of $v$. In particular, we make the left-descendant of
$v$ inherit all the $r$-close leaves by setting
$$\vec{x}' = (x_{l}, \ldots, x_{r}, \underbrace{0, \ldots, 0}_{k-r+l})$$
and let the right-descendants evenly split the remaining ones by setting
$$\vec{x}'' = (\underbrace{0, \ldots, 0,}_{r-l+1} \frac{x_{r + 1}}{2^{l} - 1},
\ldots, \frac{x_{k}}{2^{l} - 1}, \underbrace{0, \ldots, 0}_{l}).$$

To make a piecewise split useful we have to show that whenever $\vec x'$ and
$\vec x''$ are constructible, so is $\vec x$. This follows from iterated
application of Observation~\ref{obse:splitting} if the intermediate vectors
$\vec x^*$ (the leaf-vectors assigned to the intermediate vertices in the
binary tree of depth $l$) satisfy the requirement $|\vec x^*|\le d$. This
condition will be satisfied in the cases where we apply piecewise split.

\medskip

The advantage of a piecewise split is that since the coordinates with a small
index are fully given to the left-descendant, their weight is multiplied by
$2^l$. We will set the parameters such that this makes the weight of $\vec x'$
reach $1$, ensuring its constructibility by Lemma~\ref{lemma:payofflemma}. For
the right-descendants the weight-gain on the non-zero coordinates of
the assigned leaf-vector is uniformly distributed, but
tiny: only a factor $1+\frac{1}{2^l -1}$.
Furthermore the leaf-vector starts with many zeros, so we
can perform a large number of fair splits and hope that the resulting
leaf-vector is ``better'' in some way than $\vec x$, for example its weight
increases.
This will not always be the case in reality, because the behaviour of
the weight in the optimal process is more subtle and can oscillate.
This represents yet another, more serious technicality to handle in
Section~\ref{sec:FormalProof}.

The cut subroutine in the next paragraph describes more formally the above
combination of the piecewise
split and the fair splits on the right-descendants.

\paragraph{The Cut Subroutine.}
The cut subroutine has a single parameter $l$ with $1\le l\le k/2$. It can be
applied to a leaf-vector $\vec x$ with $x_i=0$ for $i<l$
and of weight
$w(x)\ge2^{-l}$. It consists of a piecewise split with parameters $l$ and $r=r(l,\vec{x})$,
where $r\ge2l-1$ is the
smallest index such that $\sum_{i=l}^{r} x_{i}/2^i \ge2^{-l}$. The choice of
$r$ ensures that the leaf-vector $\vec x'$ of the left-descendant is
constructible by Lemma~\ref{lemma:payofflemma}. Then we apply an
$(r-2l+1)$-times
iterated fair split to the leaf-vector $\vec x''$ of the right-descendants to
obtain a leaf-vector
$$C^l(\vec
x)=\left(\underbrace{0,\ldots,0}_{l},\frac{x_{r+1}}{2^{r-l}(2^l-1)},\ldots,\frac{x_k}{2^{r-l}(2^l-1)},\underbrace{0,\ldots,0}_{l},\frac{d}{2^{r-2l+1}},\ldots,\frac{d}{2}\right).$$

As we ensured the constructibility of $\vec x'$ and because we observed
similar implications for the piecewise split and fair split
operations we have that the constructibility of $C^l(\vec x)$ implies the
constructibility of $\vec x$.

Note that we maintain enough zeros at the left end of
$C^l(\vec x)$ to be able to use the cut subroutine recursively. As an illustration we do this first in the
simplest case $l=1$. This takes up the rest of this subsection and already gives a factor $2$ improvement over
the bound on $d$ which we obtained in
Section~\ref{sec:FormalDef} using only the fair split. In the next subsection we analyze what the recursive application of the cut subroutine gives us if $l$ is set to be a large constant. This is enough to obtain asymptotically tight results, but later, in Section~\ref{sec:FormalProof}, where we give a formal discrete version of our argument we chose $l$ to be logarithmic in $k$ to obtain optimal results. The discussion in the rest of this subsection is not crucial for the next subsection and the whole informal argument in this section is not necessary for the formal treatment in the next section.

We saw earlier using repeated applications of the
fair split and Observation~\ref{obse:constructibility} that in order to prove the existence of
$(k,d)$-trees
it is enough to see that the
vector $\vec x^0=(0,\ldots,0,1,2,4, \ldots d/2)$ is constructible. Our plan is
to establish this through the repeated application of cuts with parameter
$l=1$: we recursively define $\vec x^i=C^1(\vec x^{i-1})$ and hope that we
eventually obtain a vector with $w(\vec x^i)\geq 1$. By
Lemma~\ref{lemma:payofflemma} this would establish the constructibility of
$\vec x^i$, and through that, also the constructibility of $\vec x^0$ and the
existence of $(k,d)$-trees.

In order to just get started with the first cut we need
$w(\vec x^0)=\log d\cdot \frac{d}{2^{k+1}}\ge2^{-l}=1/2$
and thus $d>2^k/k$. It turns out that if $d$ is chosen slightly larger at $d=(1+\epsilon)2^k/k$
and $k$ is large enough for $\epsilon$, our plan of repeated cut
subroutines with parameter $l=1$ can indeed be carried
through.
Note that this bound on $d$ is a factor of $2$ smaller than the bound
obtained from fair split and Lemma~\ref{lemma:payofflemma} alone in the
previous section, but it is still larger by a factor of $e/2$ than the bound
we need to prove Theorem~\ref{kd-trees}. For the stronger bound we will need
cuts with parameter $l>1$, but as an illustration we give a very rough sketch
here why $l=1$ is enough if $d=(1+\epsilon)2^k/k$.

Let us start with examining the first cut producing $\vec x^1$. Except the first $\approx \log k$ coordinates, each coordinate contributes the
same $\frac{d}{2^{k+1}} = (1+\epsilon)/(2k)$ to the weight of $\vec
x^0$. Thus the first piecewise split will have parameter $r_1 \approx
k/(1+\epsilon)\approx(1-\epsilon)k$.
After the piecewise split the leaf-vector
of the right descendant will have only $\approx\epsilon k$ non-zero entries,
but the contribution to the weight of each of these entries is doubled to
$\approx(1+\epsilon)/k$.
This contribution is not changing during the repeated
fair splits, but $r_1-1$ new non-zero entries show up, each with the
``standard'' $\approx(1+\epsilon)/(2k)$ contribution to the weight. In total
we will have $w(\vec x^1)\approx(1+2\epsilon)/2$, a noticeable
$\frac{\epsilon}{2}$ improvement over
$w(\vec x^0)$.

This improvement in the weight of the coordinates towards the beginning
of the leaf-vector makes the parameter of the second cut slightly
smaller at $r_2\approx(1-2\epsilon)k$, further improving the weight of $\vec
x^2$. In general we will have $r_i\approx(1-i\epsilon)k$ and $w(\vec
x^i)\approx(1+(i+1)\epsilon)/2$. This works till $r_i>k/2$. After that
threshold (at around $i=1/(2\epsilon)$) the rate by which the weight increases
slows a little, but we will still have an index $i<2/\epsilon$ with $w(\vec
x^i)>1$ as needed to finish this argument.

\subsection{Passing to continuous} \label{subsection:cont}

The goal of this subsection is to give the continuous
motivation behind the formal discrete proof of the next section as well as
to elucidate {\em why} the discrete construction should work.
The continuous function $F$, defined in this subsection via a differential equation, is also helpful in studying the {\em size} of the constructed formulas
(see Section~\ref{sec:Outlook}).

Recall that our goal is to obtain a $(k,d)$-tree for
\begin{equation} \label{eq:defiofd}
 d = \frac{1}{T} \frac{2^{k+1}}{k}
\end{equation}
where $T$ should be as large as possible and $k$ suitably large for a given
$T$. To establish the upper bound of Theorem~\ref{kd-trees} we need this for
any $T<e$. A similar result for a value $T>e$ would contradict to the lower
bound in the same theorem. As in the analysis in the previous section we
build the $(k,d)$-tree with repeated applications of the cut subroutine. We choose
the parameter $l$ of the subroutine to be a large constant, i.e., it depends on $T$, but not on $k$.

After fixing a target constant $T$, it will be helpful to consider
the leaf-vectors in a normalized form, which then will enable us
to interpret them as continuous functions on the $[0,1]$ interval.

First we normalize the leaf-vector $\vec z=(z_0,\ldots,z_k)$ to get
$(y_0,\ldots,y_k)$ with $y_i=2^{k+1-i}z_i/d$. Note that we have $w(\vec
z)=\sum_iy_i/(kT)$, so, in particular, $\vec z$ is constructible whenever
$\sum_iy_i\geq kT$. We associate a real function $f:[0,1]\to\Real$ to the
leaf-vector $\vec z$ by the formula $f(x)=y_{\lfloor kx\rfloor}$. Although this
is a step function we will treat it as a continuous function. Clearly, as $k$
increases, this is more and more justified.

Our Lemma~\ref{lemma:payofflemma} translates to our new setting as follows.

\begin{lemma} \label{lemma:payoffContinuous}
If the real function $f$ associated to the leaf-vector $\vec z$
satisfies
$$\int_{x = 0}^{1} f(x) dx \geq T$$ then $\vec{z}$ is constructible.
\end{lemma}

We now illustrate how the cut subroutine (with sufficiently high parameter $l$
chosen for $\epsilon>0$) can be applied to achieve the target $T = e -
\epsilon$. This is an informal analysis (as everything in this
section) and hence we will allow ourselves to ignore what happens to
the function $f(x)$ on $o(1)$-long
subintervals of $[0,1]$.

Our first job is to see how a single application of the cut subroutine with
parameter $l$ transforms the real function $f$ associated to a leaf-vector. The
cut starts with a piecewise split with parameters $l$ and $r$ where $r$ is
chosen minimal with respect to the condition that the leaf-vector of the
left-descendant has weight at least $1$. Let us set $v=r/k$ and consider the
function

$$
f_{\rm left}(x) = \left\{ \begin{array}{ll} 2^lf(x) & x\in [0,v) \\
                                   0 & x\in [v,1]\end{array}\right.$$

Notice that the true transformation of the leaf-vectors involves a left shift
of $l$ places for the entries up to $x_r$ followed by zeros. This left shift
explains why the normalized vector and thus the values of the associated real
function are multiplied by $2^l$. By the same left shift the values of the real
function should also be shifted to the left by $l/k$, but as $l/k=o(1)$ this
small effect is ignored here.

Notice that by Lemma~\ref{lemma:payoffContinuous} the choice of $r$ translates to a choice of $v$ that makes
\begin{equation} \label{eq:cont-definition}
T = \int_0^1f_{\rm left}(x)\enspace dx = 2^l \int_0^vf(x)\enspace dx.
\end{equation}

The real function $f_{\rm right}$ associated to the  leaf-vector of the right-descendants of the
piecewise split with parameters $l$ and $r$ can similarly be approximated by

$$f_{\rm right}(x)= \left\{ \begin{array}{ll}0&x\in[0,v)\\
\frac{2^l}{2^l-1}f(x) & x\in [v,1]\end{array}\right.$$

Recall that the cut subroutine further applies the $(r-l)$-times iterated fair split
to the leaf-vector. Because of our normalization fair splits do not change the
values of the associated real function, but they move these values to the
left in the domain of the function. In our case the amount of this translation
is $(r-l)/k\approx v$. The fair splits also introduce constant $1$ values in
the right end of the $[0,1]$ interval freed up by the shift. This explains why
the following function approximates well the resulting real function:
\begin{equation}
\label{eq:right-cont}
C^l_f(x)=\left\{ \begin{array}{ll} \frac{2^l}{2^l-1}f(x+v) & x\in [0,1-v) \\
                                   1 & x\in [1-v, 1]
\end{array} \right.
\end{equation}

In the following we analyze how the repeated application of these cut
subroutines changes the associated real function. For this analysis
we define a two-variable function $F(t,x)$ that approximates well how the
function develops. Here $t\ge0$ represents the ``time'' that has
elapsed since we started our process, that is the cumulative length of the
shifts $v$ we made. The value $x\in [0,1]$ stands for the single variable of
our current real function.
In other words, for each fixed $t$, $F(t,x)$ should be a good approximation of the real function associated to the leaf-vector after
$t/v_{\rm avg}$ infinitesimally small cuts were made
(where $v_{\rm avg}$ is the length of the average cut).
We have the initial
condition $F(0,x)=1$ as the constant $1$ function approximates the real
function associated to our original leaf-vector of
$(0,\ldots,0,1,2,\ldots,d/2)$ (in fact, the two functions coincide except for
the $o(1)$ length subinterval of their domain where the latter function is
$0$). We also have $F(t,1) = 1$ for every $t\geq 0$ by \eqref{eq:right-cont}.

We make yet another simplification: we assume that $l$ is chosen large
enough making $v$ so small, that we can treat it as infinitesimal.
We will denote the
cut parameter $v$ at time $t$ by $v_t$. Recall that
by~\eqref{eq:cont-definition} $v_t$ can be defined by the formula
$\int_0^{v_t}F(t,x)\enspace dx=T/2^l$, which simplifies to
$v_t\approx T/(2^lF(t,0))$ after approximating the integral.

For the ''right-descendant'' we have by \eqref{eq:right-cont} that
$$F(t+v_t,x)=C^l_{F(t,\cdot)}(x)=\frac{2^l}{2^l-1}F(t,x+v_t)$$
if $x<1-v_t$ and $F(t+v_t,x)=1$ for every $1-v_t\leq x \leq 1$.
Using the approximation $\frac{2^l}{2^l-1}\approx1+\frac1{2^l}$ (justified as
$l$ is considered ``large'') we approximate
the above equation with
$$F(t+v_t,x)\approx   \left(1 + \frac{1}{2^{l}} \right) F(t,x+v_t) \approx \left(1+ \frac{v_t F(t,0)}{T}\right)F(t,x+v_t).$$
This gives us an equation on the derivative of $F(t,x)$ in direction $(1,-1)$.
For any $s>1$, define the function $F_s(t):[s-1,s]\rightarrow {\mathbb R}$ by $F_s(t) = F(t, s-t)$.
Then rewriting the above with $s=x+t+v_t$, we have
$$F_s(t+v_t) \approx \left(1+ \frac{v_tF(t,0)}{T}\right)F_s(t).$$
\begin{equation}\label{eq:derivative}
\frac{F'_s(t)}{F_s(t)} \approx \frac{F_s(t+v_t) - F_s(t)}{v_tF_s(t)} \approx
\frac{F(t,0)}{T}.
\end{equation}
Integrating  we obtain
$$ \int_{s-1}^s (\ln F_s(t))'dt \approx \frac{1}{T} \int_{s-1}^s F(t,0) dt.$$
The left hand side evaluates to $\ln F_s(s)=\ln F (s,0)$ by the boundary
condition $F_s(s-1) = 1$.

As our last simplifying assumption we assume that the function $F(t,0)$
increases monotonically. Therefore the right hand side is at least
$F(s - 1,0)/T$, which implies that
\begin{equation} \label{eq:expressingfunctiononsbypreviousfunctiondiverges}
F(s,0) \geq e^{\frac{F(s-1,0)}{T}}.
\end{equation}

Now the increasing function $F(t,0)$ either goes to infinity or tends to a
finite limit $a$. If the latter happens we have
$a\geq e^{a/T}$.
But classic calculus shows that $a\le e^{a/e}$ for all real
$a$, so this implies $T\ge e$. In the case $T<e$, that we study here, $F(t,0)$
must then tend to infinity. From the assumed monotonicity of $F(t,0)$ and
\eqref{eq:derivative} we obtain $F(t,x) > \sqrt{F(t,0)}$ if $1/2\leq x \leq 1$.
Since $F(t,0)$ tends to infinity, so does
$\int_0^1F(t,x)\enspace dx$. Whenever this integral grows above $T$,
the current leaf-vector is constructible by
Lemma~\ref{lemma:payoffContinuous} finishing our highly informal proof of the
existence of $(k,d)$-trees.

The above continuous heuristic is the underlying idea of
the construction described in the next section.
It provides a good approximation to what happens in the discrete case.
Instead of dealing with all the introduced approximation errors in a
precise manner, we give a direct discretized proof where
we explicitly make sure
that our many simplifying assumptions are satisfied and the
approximations are correct.


\section{Formal Construction of $(k,d)$-Trees} \label{sec:FormalProof}

In this section we complete the proof of Theorem~\ref{kd-trees}.

Before proving Theorem~\ref{kd-trees} we first set up two of the main
ingredients of our construction.
Let us fix the positive integers $k$, $d$ and $l$. To simplify the
notation, we will not show the
dependence on these parameters in the next definitions, although $d'$,
$E$, $C_r$ and $C^*_r$ all depend on them. We let
$$d'=d\left(1-\frac{1}{2^l-1}\right).$$
For a vector $\vec x=(x_0,\ldots,x_k)$ we define
$$E(\vec x)=\left(\lfloor x_1/2\rfloor,\lfloor x_2/2\rfloor,\ldots,\lfloor
x_k/2\rfloor,\lfloor d'/2\rfloor\right).$$
We denote by $E^m(\vec x)$ the vector obtained
from $\vec x$ by $m$ applications of the operation $E$.
Using the simple observation that
$\lfloor\lfloor a\rfloor/j\rfloor=\lfloor a/j\rfloor$ if $a$ is real and $j$
is a positive integer we can ignore all roundings but the last.
$$E^m(\vec x)=\left( \left\lfloor \frac{x_m}{2^m}\right\rfloor, \left\lfloor
    \frac{x_{m+1}}{2^m}\right\rfloor, \ldots, \left\lfloor
\frac{x_k}{2^m}\right\rfloor, \left\lfloor \frac{d'}{2^m}\right\rfloor,
\left\lfloor \frac{d'}{2^{m-1}}\right\rfloor, \ldots, \left\lfloor
  \frac{d'}{2}\right\rfloor\right).$$

For $l\le r\le k$ and the vector $\vec x$ as above we define the $(k+1)$-tuples
$C_r(\vec x)$ and $C^*_r(\vec x)$ by the following formulas:
\begin{align*}
C_r(\vec x) &= \Big( \underbrace{0,\ldots ,0,}_{r+1-l}
\underbrace{\left\lfloor\frac{x_{r+1}}{2^{l}-1}\right\rfloor, \ldots ,
  \left\lfloor\frac{x_{k}}{2^{l}-1}\right\rfloor,}_{k-r}
\underbrace{\left\lfloor \frac{d'}{2^l}\right\rfloor, \ldots, \left\lfloor
    \frac{d'}{2}\right\rfloor}_{l} \Big)\\
C^*_r(\vec x) & = ( \underbrace{x_l, x_{l+1},\ldots, x_r,}_{r+1-l} \underbrace{0,0,
  \ldots , 0}_{k-r+l} ).
\end{align*}

Note that for the following lemma to hold we could use $d$ instead of $d'$ in
the definition of $E$, and we could also raise most of the constant terms in
the definition of $C_r$. The
one term we cannot raise is the entry $\lfloor d'/2^l\rfloor$ of
$C_r(\vec x)$ right after $\lfloor x_k/(2^l-1)\rfloor$. If we used a higher value
there, then one of the children of the root of the tree constructed in the
proof below would have more than $d$ leaves $k$-close. We use $d'$
everywhere to be consistent and provide for
the monotonicity necessary in the proof of Theorem~\ref{kd-trees}.

The first part of the next lemma states the properties of our  fair
split procedure
(which is somewhat modified compared to the informal treatment of Section~\ref{sec:Informal}),
the second part does the same for the cut
subroutine.

\begin{lemma}\label{split}
Let $k$, $d$ and $l$ be positive integers and $\vec x\in \N^{k+1}$ with $|\vec x|\le
d$.
\begin{enumerate}[(a)]
\item$|E(\vec x)|\le d$. If $E(\vec x)$ is $(k,d)$-constructible, then so is $\vec x$.
\item For $l\le r\le k$ we have $|C_r(\vec x)|\le d$ and $|C^*_r(\vec x)|\le d$. If both
  of $C_r(\vec x)$ and $C^*_r(\vec x)$ are $(k,d)$-constructible and
$|C^*_r(\vec x)|\le d/2^l$, then $x$ is also $(k,d)$-constructible.
\end{enumerate}
\end{lemma}

\begin{proof} (a) We have $|E(\vec x)|\le|\vec x|/2+d'/2<d$. If there exists a
$(k,d,E(\vec x))$-tree, take two disjoint copies
of such a tree and connect them with a new root vertex, whose children are the
roots of these trees. The new binary tree so obtained is a $(k,d,\vec x)$-tree.

(b) The sum of the first $k+1-l$ entries of $C_r(\vec x)$ is at most
$|\vec x|/(2^l-1)\le d/(2^l-1)$, the remaining fixed terms sum to less than
$d'=d(1-1/(2^l-1))$, so $|C_r(\vec x)|\le d$. We trivially have
$|C^*_r(\vec x)|\le|\vec x|\le d$.

Let $T$ be a $(k,d,C_r(\vec x))$-tree and $T^*$ a $(k,d,C^*_r(\vec x))$-tree. Consider a
full binary tree of height $l$ and attach $T^*$ to one of the $2^l$
leaves of this tree and attach a separate copy of $T$ to all remaining
$2^l-1$ leaves. This way we obtain a finite binary tree $T'$. We claim that
$T'$ is a $(k,d,\vec x)$-tree showing the constructibility of $\vec x$ and finishing the
proof of the lemma. To check condition (i) of the
definition of a $(k,d,\vec x)$-tree notice that no leaf of $T'$ is in distance less
than $l$ from the
root, leaves in distance $l\le j\le r$ are all in $T^*$ and leaves in distance
$r<j\le k$ are all in the $2^l-1$ copies of $T$. Hence
$$\left(0,\ldots , 0, x_l, \ldots , x_r,  (2^l
  -1)\left\lfloor\frac{x_{r+1}}{2^{l}-1}\right\rfloor, \ldots ,
  (2^l-1)\left\lfloor\frac{x_{k}}{2^{l}-1}\right\rfloor\right)$$
is a leaf-vector of the root of $T'$, thus $\vec x$ is also a leaf-vector.

Condition (ii) is satisfied for the root, because it has at most
$|\vec x|\leq d$ leaves that are $k$-close.
Notice that the nodes of $T'$ of distance at least $l$ from the root are also
nodes of $T^*$ or a copy of $T$, so they satisfy condition (ii). There are
two types of vertices in distance $0<j<l$ from the root. One of them has
$2^{l-j}$ copies of $T$
below it, the other has one less and also $T^*$.
In the first case we can bound the number of $k$-close leaves by
\begin{eqnarray}
& &2^{l-j} \left( \frac{x_{r + 1} + \ldots + x_{k}}{2^{l} - 1} +
  \frac{d'}{2^{l}} + \ldots + \frac{d'}{2^{l-j+1}} \right)
\leq   \frac{2^{l-j}}{2^{l} - 1} \cdot d + d' \left(1 -
  \frac{1}{2^{j}}\right) \nonumber \\
& = & \frac{d}{2^{l} - 1} \left(2^{l - j} + (2^{l} -2) \left(1 -
    \frac{1}{2^{j}}\right)  \right)
\leq \frac{d}{2^{l} - 1} (2^{l} - 1) \nonumber \\
& = & d, \label{eq:noteexpression}
\end{eqnarray}
where the last inequality holds since $j\geq 1$.
We also point out that here it is crucial that we use the specific value of $d'<d$.

In the second case, the number of $k$-close leaves is bounded by
\begin{eqnarray*}
& & (2^{l-j} - 1) \left( \frac{x_{r + 1} + \ldots + x_{k}}{2^{l} - 1} + \frac{d'}{2^{l}} + \ldots + \frac{d'}{2^{l-j+1}} \right) + |C^*_r(\vec x)| \\
& \leq & (2^{l-j} -1) \frac{d}{2^{l-j}} + \frac{d}{2^{l}} \leq d,
\end{eqnarray*}
where we used the result of the calculation in \eqref{eq:noteexpression}.
\end{proof}

Armed with the last lemma we are ready to prove Theorem~\ref{kd-trees}.
\begin{proof}[Proof of Theorem~\ref{kd-trees}]
We will show that
$(k,d)$-trees exist for large enough $k$ and
$d=\lfloor2^{k+1}/(ek)+100\cdot2^{k+1}/k^{3/2}\rfloor$.

We set $l=\lfloor\log (k)/2\rfloor$,  so $2^l\sim\sqrt k$. 
We will construct our vectors in such a way that they start with $2l + 1$ zeros.

We define the vectors $\vec x^{(t)}=(x^{(t)}_0,\ldots,x^{(t)}_k)\in\N^{k+1}$
recursively. We start with $\vec x^{(0)}=E^{k-2l}(\vec z)$, where $\vec z$ denotes the
vector consisting of $k+1$ zeros. For $t\ge0$ we define
$\vec x^{(t+1)}=E^{r_t-3l}(C_{r_t}(\vec x^{(t)}))$, where $r_t$ is the smallest
index in the range $3l\le r_t\le k$ with
$\sum_{j=0}^{r_t}x^{(t)}_j/2^j\ge2^{-l}$. At this point we may consider the
sequence of the vectors $\vec x^{(t)}$ end whenever the weight of one of
them falls below $2^{-l}$ and thus the definition of $r_t$ does not make
sense. But we will establish below that this never happens and the sequence is
infinite.

Notice first, that for all $t$, the first $2l+1$ entries of $\vec x^{(t)}$ are zeros and all other entries are obtained from $d'$ by repeated
application of dividing by an integer (namely by $2^l-1$ or by a power of $2$)
and taking lower integer part. As we have observed earlier in this section, we
can ignore all roundings but the last. This way we can write each of these
entries in the form
$\left\lfloor\frac{d'}{2^i(2^l-1)^c}\right\rfloor=\left\lfloor\frac{d'}{2^{i+lc}}\alpha^c\right\rfloor$ for some
non-negative integers $i$ and $c$ and with $\alpha=1+1/(2^l-1)$. 
We will use the values $q_t=r_t-2l$ (the length of the ``left shift'' when going from 
$\vec x^{(t)}$ to $\vec x^{(t+1)}$) to give the
exponents of $2$ and $\alpha$ explicitly for all $2l<j\le k$ and all $t$.  
We claim that
\begin{equation}\label{numeric}
x_j^{(t)}=\left\lfloor\frac{d'}{2^{k+1-j}}\alpha^{c(t,j)}\right\rfloor
\end{equation}
where $c(t,j)$ is the largest integer $0\le
c\le t$ satisfying $\sum_{i=t-c}^{t-1}q_i\le k-j$, $c(0,j)=0$ for
all $j$, and $c(t,j) = 0$ if $q_{t-1} > k-j$.

The formal inductive proof of this formula becomes straightforward 
once we adopt the right perspective to keep track of
how the entries develop during the evolution of the 
vectors from $\vec x^{(0)}$ to $\vec x^{(t)}$. Namely, that each current entry
$x_j^{(t)}$ has an ``ancestor'' that entered close to the right end of some vector 
$\vec x^{(t')}$ as $x^{(t')}_{j'}=\lfloor d'/2^{k+1-j'}\rfloor$ and then got shifted to its current position as a result of 
applications of a certain number of $E$s and some $C_{r_i}$s.
As a first approximation we count a division by $2$ for each time the entry moved 
one place to the left --- this explains the exponent $k+1-j$ of the $2$ in
the formula. This is an accurate assumption whenever the left-shift was the result of an
application of $E$. However when some $C_{r_i}$ was applied to an already
existing entry, then the entry moved $l$ 
places to the left and the actual division was by $2^l-1$ instead of
$2^l$, so an extra factor of $\alpha$ has to be included for correction.
The exponent $c(t,j)$ counts exactly how many such $\alpha$-factors are accumulated:
this is exactly $t-t'$, that is, the number of $C_{r_i}$ that were applied to the ancestor of $x_j^{(t)}$ after  it got introduced into $\vec x^{(t')}$. 
Starting from the rightmost entry, the entry must be shifted $(k-j)$ places to the left 
and these  are accumulated via the left-shifts $q_{t-1}, \ldots, q_{t'}$, that are the
result of the applications of $E^{r_{t-1} -3l} \circ C_{r_{t-1}}, \ldots E^{r_t' -3l}\circ C_{r_{t'}}$, respectively.

We claim next that $c(t,j)$ and $x^{(t)}_j$ increase monotonously in $t$ for
each fixed $2l<j\le k$, while $q_t$ decreases monotonously in $t$. We prove
these statements by induction on $t$. We have
$c(0,j)=0$ for all $j$, so $c(1,j)\ge c(0,j)$. If $c(t+1,j)\ge
c(t,j)$ for all $j$, then all entries of $\vec x^{(t+1)}$ dominate the
corresponding entries of $x^{(t)}$ by \eqref{numeric}. If
$x^{(t+1)}_j\ge x^{(t)}_j$ for all $j$, then we have $r_{t+1}\le r_t$ by the
definition of these numbers, so we also have $q_{t+1}\le q_t$. Finally, if
$q_0 \geq q_1 \geq \cdots \geq q_{t+1}$, then by the definition of $c(i,j)$ we
have $c(t+2,j)\ge c(t+1,j)$.

The monotonicity just established also implies that the weight of $\vec x^{(t)}$ is
also increasing, so if the weight of $\vec x^{(0)}$ is at least $2^{-l}$, then so
is the weight of all the other $\vec x^{(t)}$, and thus the sequence is infinite.
The weight of $\vec x^{(0)}$ is
\begin{align*}
\sum_{j=2l+1}^k\frac{\left\lfloor\frac{d'}{2^{k+1-j}}\right\rfloor}{2^j}
&>\sum_{j=2l+1}^k\frac{\frac{d'}{2^{k+1-j}}-1}{2^j}\\
&>(k-2l)\frac{d'}{2^{k+1}}-2^{-2l}\\
&>(k- \log k)\frac{1-\frac1{2^{l}-1}}{ek}-2^{-\log k +2} \ \longrightarrow \ \frac{1}{e},
\end{align*}
where the last inequality follows from $d>2^{k+1}/(ek)$ and the last term
tends to $e^{-1}$ as $k$ tends to infinity, so it is larger than $2^{-l}$ for
large enough $k$.

We have just established that the sequence $\vec x^{(t)}$ is
infinite and coordinate-wise increasing. Notice also that these vectors were
obtained through the operations $E$ and $C_r$ from the all-zero vector, so by
Lemma~\ref{split} we must have $|\vec x^{(t)}|\le d$ for all $t$. Therefore
the sequence $\vec x^{(t)}$ must stabilize eventually. In fact, as the sequence
stabilizes as soon as two consecutive vectors agree, it must stabilize
in at most $d$ steps. So for some fixed vector $\vec x=(x_0,\ldots,x_k)$ we have
$\vec x^{(t)}=\vec x$ for all $t\ge d$. This
implies that $q_t$ also stabilizes with $q_t=q$ for $t\ge
d$. Equation \eqref{numeric} as applied to $t>d+k$ simplifies to
\begin{equation}\label{num}
x_j=\left\lfloor\frac{d'}{2^{k+1-j}}\alpha^{\lfloor(k-j)/q\rfloor}\right\rfloor.
\end{equation}

Recall that $q=q_t=r_t-2l$ (take $t\ge d$), and $r_t$ is defined
as the smallest index in the range $3l\le r\le k$ with
$\sum_{j=0}^rx_j/2^j\ge2^{-l}$. Thus we have $q\ge l$.
\begin{prop}\label{ql}
We have $q = l$.
\end{prop}
\begin{proof} Assume for contradiction that $q>l$. Then by the minimality of $r_t$
we must have
\begin{eqnarray*}
2^{-l}&>&\sum_{j=0}^{2l+q-1}\frac{x_j}{2^j}\\
&=&\sum_{j=2l+1}^{2l+q-1}\frac{\left\lfloor\frac{d'}{2^{k+1-j}}\alpha^{\lfloor(k-j)/q\rfloor}\right\rfloor}{2^j}\\
&>&\sum_{j=2l+1}^{2l+q-1}\frac{\frac{d'}{2^{k+1-j}}\alpha^{(k-j)/q-1}-1}{2^j}\\
&>&(q-1)\frac{d'}{2^{k+1}}\alpha^{k/q-4}-2^{-2l}.
\end{eqnarray*}
In the last inequality we used $\frac{j}{q} \leq \frac{2l+q}{q} = 1 + \frac{2l}{q} \leq 1 + \frac{2l}{l} = 3$. This inequality simplifies to
\begin{equation}\label{impossible}
2^{-l}(1+2^{-l})\alpha^4\frac{2^{k+1}}{d'}>(q-1)\alpha^{k/q}.
\end{equation}
We consider the problem of minimizing the right hand side over {\sl real} numbers $q > 2$:
Simple calculus gives that the corresponding value of $q$
is sandwiched 
between $q^*-1$ and $q^*-2$, where $q^*=q^*(k)=k\ln\alpha \sim \sqrt{k}$, and this minimum is more than
$(q^*-3)e$. Using $\alpha=1+1/(2^l-1)>e^{2^{-l}}$ we have $q^*\ge k/2^l$.
So \eqref{impossible} yields
$$2^{-l}(1+2^{-l})\alpha^4\frac{2^{k+1}}{d'}>\frac{ke}{2^l}-3e.$$
Thus,
\begin{displaymath}
(1+2^{-l})\alpha^4\frac{2^{k+1}}{d'} - ke + 3e2^{l} > 0.
\end{displaymath}
Substituting our choice for $d$, $d'$, $\alpha$ and $l$ (as functions of $k$)
and assuming that $k$ is large enough we get that the left hand side is at most
{\allowdisplaybreaks
\begin{eqnarray*}
& & \left(1 + \frac{2}{\sqrt{k}} \right) \left(1 + \frac{4}{\sqrt{k}} \right)^{4} \frac{2^{k+1}}{d\left(1 - \frac{4}{\sqrt{k}}\right)} - ek + 3e\sqrt{k} \\
& \leq & \left(1 + \frac{19}{\sqrt{k}} \right) \frac{ek}{\left(1 - \frac{4}{\sqrt{k}}\right)\left(1 + \frac{99e}{\sqrt{k}}\right)} - ek + 3e\sqrt{k} \\
& \leq & \left(1 - \frac{200}{\sqrt{k}}\right)ek - ek + 3e\sqrt{k}  =  - 200e\sqrt{k} + 3e\sqrt{k} < 0,
\end{eqnarray*}}
which is a contradiction.
Hence $q=l$ as claimed.
\end{proof}

Having proved Proposition~\ref{ql} we return to the proof of Theorem~\ref{kd-trees}.
We establish by downward induction that the vectors
$\vec x^{(t)}$ are constructible. We start with
$t=d$.
Using \eqref{num} and the fact that $q = l$ (Proposition~\ref{ql}) we get that
for large $k$
{\allowdisplaybreaks
\begin{eqnarray*}
 w(\vec x^{(d)}) & \geq & \frac{x^{(d)}_{2l + 1}}{2^{2l + 1}} \geq \frac{d'}{2^{k + 1}} \alpha^{\frac{k - 2l - 1}{l} - 1} - 1 \\
      & \geq & \frac{d'}{2^{k + 1}} \alpha^{\frac{k}{2l}} - 1 \\
      & \geq & \frac{d/2}{2^{k + 1}} \left(1 + \frac{1}{\sqrt{k}} \right)^{\frac{k}{\log k}} - 1 \\
      & \geq & \frac{1}{2ek} e^{\frac{k}{2 \sqrt{k} \log k}} - 1 > 1.
\end{eqnarray*}}
So by Lemma~\ref{lemma:payofflemma}, $\vec x^{(d)}$ is constructible.

Now assume that $\vec x^{(t+1)}$ is constructible for some $t\leq d-1$.
Recall that $\vec x^{(t+1)}=E^{r_t-3l}(C_{r_t}(\vec x^{(t)}))$, so by (the
repeated use of) part (a) of Lemma~\ref{split}, $C_{r_t}(\vec x^{(t)})$ is
constructible. By part (b) of the same lemma, $\vec x^{(t)}$ is also
constructible (and thus the inductive step is complete) if we can (i) show
that $C^*_{r_t}(\vec x^{(t)})$ is constructible and (ii) establish that
$|C^*_{r_t}(\vec x^{(t)})|\le d/2^l$. For (i) we use the definition of $r_t$:
$\sum_{j=0}^{r_t}x^{(t)}_j/2^j\ge2^{-l}$. But the weight of
$C^*_{r_t}(\vec x^{(t)})$ is $\sum_{j=l}^{r_t}x^{(t)}_j/2^{j-l}$, so the
contribution of each term with $j\ge l$ is multiplied by $2^l$, while the
missing terms $j<l$ contributed zero anyway as the first $2l+1 > l$
coordinates of $x^{(t)}$ are $0$. This shows that the
weight of
$C^*_{r_t}(\vec x^{(t)})$ is at least $1$ and therefore Lemma~\ref{lemma:payofflemma} proves
(i).

For
(ii) we use monotonicity to see $|C^*_{r_t}(\vec x^{(t)})|=\sum_{j=2l+1}^{r_t}x^{(t)}_j\le\sum_{j=2l+1}^{r_t}x_j$.
Here $r_t\le r_0\le k/2$ for large enough $k$. Using the
fact that $q=l$ (Proposition~\ref{ql}) and \eqref{num} we further have for
large $k$ that
{\allowdisplaybreaks
\begin{eqnarray*}
 |C^*_{r_t}(\vec x^{(t)})| & \leq & d' \sum_{j = 2l + 1}^{r_{t}} \frac{1}{2^{k + 1 - j}} \cdot \alpha^{\frac{k-j}{l}} \\
         &  \leq & d \sum_{j = 2l + 1}^{r_{t}} \left( \frac{\alpha}{2} \right)^{k-j} \\
	 &  \leq & d \left(\frac{3}{4}\right)^{k - r_t} \cdot 4 \enspace \enspace \left(\text{since $\frac{\alpha}{2} \leq \frac{3}{4}$}\right) \\
         &  \leq & d \left(\frac{3}{4}\right)^{\frac{k}{2}} \cdot 4 \enspace \enspace \left(\text{since $r_{t} \leq \frac{k}{2}$}\right) \\
	 &  < & \frac{d}{\sqrt{k}} \leq \frac{d}{2^{l}}.
\end{eqnarray*}}
This finishes the proof of (ii) and hence the
inductive proof that $\vec x^{(t)}$ is constructible for every $t$.

As $\vec x^{(0)}=E^{k-2l}(\vec z)$ is constructible, Lemma~\ref{split}~(a) implies that
the all-zero vector $\vec z$ is also constructible. Thus the larger vector
$(0,\ldots,0,d)$ is constructible too, and by
Observation~\ref{obse:constructibility} there exists a $(k,d)$-tree.
\end{proof}

\section{Proof of the Lower Bound of Theorem \ref{main}} \label{sec:ProofOfLowerBound} \label{sec:lowerboundproof}

Let us first give the intuition of where the factor two improvement is
coming from compared to using just the classical Local Lemma (c.f. \eqref{kst})
and what is the reason behind the counterintuitive probabilities
we use to choose the assignment.

When using the classical Local Lemma to give a lower bound on $f(k)$,
one sets each variable independently to {\tt true} or {\tt false} uniformly at random
and the $\left\lfloor\frac{2^k}{ek}\right\rfloor$ lower bound in \eqref{kst} is immediate.
However in the proof we must consider all pairs of clauses sharing a variable as
``dependent on each other''.  The lopsided version of the Local Lemma
\cite{ErdosSpencer} allows for a more restricted definition of
``intersecting'' clauses. Namely, one can
consider two clauses intersect only if they contain a common variable with
{\em different sign} and this still allows for the same conclusion as in the
classical Local Lemma. If all variables in a
$(k,s)$-CNF are {\em balanced}, that is they appear an equal number of
times with either sign, then each clause intersects only at most $ks/2$
other clauses in this restricted sense, instead of the at most $k(s-1)$
other clauses it may intersect in the original sense and the factor two
improvement is immediate.
To handle the unbalanced case we consider a distribution on
assignments where the variables are assigned {\tt true} or {\tt false} values with some
bias. It would be natural to favor the assignment that satisfies more
clauses, but the opposite turns out to be the distribution that works.
This is because the clauses with many variables receiving the {\em less
frequent sign} are those that intersect more than the average number of
other clauses,
so for the use of the Lopsided Local Lemma those are the ones whose satisfiability should be boosted with the bias put on the assignments.

\begin{proof}[Proof of lower bound in Theorem \ref{main}]

Let $F$ be a $(k, s)$-CNF formula with
$s=\left\lfloor\frac{2^{k+1}}{e(k+1)}\right\rfloor$.

For a literal $v$ we denote by $d_v$ the number of occurrences of $v$ in
$F$. We set a variable $x$ to {\tt true} with probability
$P_x=\frac12+\frac{2d_{\bar x}-s}{2sk}$. This makes the negated version $\bar
x$ satisfied with probability $P_{\bar x}=\frac12-\frac{2d_{\bar
    x}-s}{2sk}\ge\frac12+\frac{2d_x-s}{2sk}$ as we have $d_x+d_{\bar x}\le
s$. So any literal $v$ is satisfied with probability at least
$\frac12+\frac{2d_{\bar v}-s}{2sk}$.

In order to prove that the formula $F$ is satisfied with non-zero probability we need
the Lopsided Local Lemma of Erd\H os and Spencer.
\begin{lemma}[\cite{ErdosSpencer}]
Let $\{ A_C\}_{C\in I}$ be a finite set of events in some probability space.
Let $\Gamma(C)$ be a subset of $I$ for each $C\in I$ such that for every
subset $J\subseteq I\setminus (\Gamma (C) \cup \{ C\})$ we have
$$Pr (A_C | \wedge_{D\in J} \bar{A}_D) \leq Pr (A_C).$$
Suppose there are real numbers $0<x_C<1$ for $C\in I$
such that for every $C\in I$ we have
$$ Pr (A_C) \leq x_C \prod_{D\in \Gamma (C)} (1-x_D).$$
Then
$$Pr (\wedge_{C\in I}\bar{A}_C) > 0.$$
\end{lemma}

For each clause $C$ in $F$, we define the ``bad event'' $A_C$
to be that $C$ is not satisfied.
Moreover,
for every $C$ in $F$ we define $\Gamma(C)$ to be the family of clauses
$D$ in $F$ that have at least one such variable in common with $C$ whose
sign is different in $C$ and $D$.
Finally, we set the value of each $x_C$ to be $x=\frac{e}{2^k}$.

We need to check that for every
subset $J\subseteq I\setminus (\Gamma (C) \cup \{ C\})$ we have
$$Pr (A_C | \wedge_{D\in J} \bar{A}_D) \leq Pr (A_C).$$
This is equivalent to $Pr (\wedge_{D\in J} \bar{A}_D | A_{C}) \leq Pr
(\wedge_{D\in J} \bar{A}_D)$. The right hand side is simply the probability of
a random assignment satisfying all clauses in $J$. The left hand side can be
interpreted as the same probability after a random assignment is modified by
setting all literals in $C$ to {\tt false}.
The random assignment does not satisfy any clauses in $J$ after the
modification that it did not satisfy before, since no clause from $J$ contains
a variable of $C$ in the opposite form.
Hence the probability of satisfying all of them does not grow by the
modification proving the inequality. (The probability might, in fact, decrease
if some of the clauses in $J$ contain a literal also present in $C$.)

We need to check also the other condition of the lemma.
Let $C$ be an arbitrary clause and let us denote the literals it contains by
$v_1,\ldots,v_k$. For $C$ not to be satisfied we must not set any of the
independent literals in $C$ to {\tt true}, and therefore we have
\begin{eqnarray*}
Pr (A_C) & = & \prod_{i=1}^k(1-P_{v_i})\\
& \le &\prod_{i=1}^k\left(\frac{1}{2}-\frac{2d_{\bar v_i}-s}{2sk}\right)\\
& \le
&\frac1{2^k}\prod_{i=1}^k\left(\left(1+\frac1k\right)\left(1-\frac{ed_{\bar
        v_i}}{2^k}\right)\right)\\
& \le & \frac{\left(1+\frac1k\right)^k}{2^k}\prod_{i=1}^k(1-x)^{d_{\bar
    v_i}}\\
& < & \frac e{2^k}(1-x)^{|\Gamma(C)|} \\
& = & x \prod_{D\in \Gamma (C)} (1-x).
\end{eqnarray*}
The inequality in the fourth line holds due to the well-known inequality that
$1 - ax \leq (1 - x)^{a}$ for every $0 < x < 1$ and $a \geq 1$.

As the conditions of the Lopsided Local Lemma are satisfied, its conclusion
must also hold. It states that the random evaluation of the variables we
consider satisfies the $(k,s)$-CNF $F$ with positive probability. Thus $F$
must be satisfiable and we have $f(k)\ge
s=\left\lfloor\frac{2^{k+1}}{e(k+1)}\right\rfloor$.

\end{proof}

\section{More About the Class MU(1) and the Neighborhood Conjecture} \label{sec:Outlook}

\subsection{Constructing binary trees and MU(1) formulas}

The structure of MU(1) formulas is well understood and it is closely related to
binary trees. (Recall that by a \emph{binary tree} we always mean a rooted tree where every non-leaf node has exactly two children.) In particular, given any binary tree $T$, we associate with
it certain CNF formulas. Similarly to Section \ref{subsection:formulas} we start with assigning distinct literals to the
vertices, assigning the negated and non-negated form of the same variable to
the two children of any non-leaf vertex.
We do not assign any literal to the
root. For each leaf of $T$ select a clause that is the disjunction of {\em some}
literals along the path from the root to that leaf
and consider the CNF formula
$F$ that is the conjunction of one clause for each leaf. Clearly, $F$ is
unsatisfiable and it has one more clause than the number of variables
associated with $T$. Note that for proving the upper bound in
Theorem~\ref{main} we used
$(k,d)$-trees $T$ and the associated unsatisfiable
(cf. Observation~\ref{obse:Funsat})
$(k,d)$-CNF $F_k(T)$ that was constructed similarly, but selecting the $k$
vertices farthest from the root on every root-leaf path.

As proved in \cite{DDK}, $F$ is a MU(1) formula if and only
if all literals associated to vertices of $T$ {\em do} appear in $F$,
furthermore
every formula in MU(1) can be obtained from a suitable binary tree this way.

Recall that $\f(k)$ denotes the smallest integer $d$ such that a $(k,d)$-tree exists. Clearly,
$f(k) \leq f_{1}(k)<f_{\tree}(k)$, and we showed that $f(k)=(1+o(1))f_{\tree}(k)$.

\subsection{On the size of unsatisfiable formulas}

By the {\em size} of a rooted tree we mean the number of its leaves and by the
{\em size} of a CNF formula we mean the number of its clauses. With this
notation the size of a containment-minimal $(k,d)$-tree $T$ and the size of
the corresponding $(k,d)$-CNF $F_k(T)$ in MU(1) are the same.

When proving the upper bound of Theorem  ~\ref{main} we constructed
$(k,d)$-trees for $d\approx\frac{2^{k+1}}{ek}$. Their size and therefore the
size of the corresponding $(k,d)$-CNF in MU(1) is at most
$2^h$, where $h$ is the height of the tree. In fact, the sizes of the trees we
constructed are very close to this
upper bound.  Therefore it makes sense to take a closer look at the height.

Recall that we associated with a vertex $v$ of a $(k,d)$-tree the vector
$(x_0,\ldots,x_k)$, where $x_j$ is the number of leaf-descendants of $v$ of
distance $j$ from
$v$. If a $(k,d)$-tree has minimal size it has no two
vertices along the same branch with identical vectors. In fact, this statement
is also true even if one forgets about the last entry in the vector. So the
minimal height of a $(k,d)$-tree is limited by the number of vectors in $\N^k$
with $L_1$ norm at most $d$, which is ${d+k\choose k}\le d^k$.
For $d=\f(k)$ this is $2^{O(k^2)}$. For the minimal size of a $(k,d)$-tree
this implies a $2^{2^{O(k^2)}}$ bound that applies whenever such a tree
exists. The same bound for minimal size
$(k,d)$-CNF formulas in MU(1) is implicit in
\cite{HSCompMU}. There is numerical evidence that the sizes of the minimal
$(k,\f(k))$-tree and the minimal $(k,f_1(k)+1)$-CNF in MU(1) might indeed be
doubly exponential in $k$ (consider the size of the minimal $(7,18)$-tree and
the minimal $(7,18)$-CNF in MU(1) mentioned below).

A closer analysis of the proof of the upper bound in Theorem~\ref{kd-trees}
shows that the height
of the $(k,d)$-tree constructed there is at most $d\log k$.
While this
is better than the general upper bound above it still allows for trees with
sizes that are doubly exponential in $k$.

This height can, however, be substantially decreased if we allow the error term
in $d$ to slightly grow. If we allow $d=(1+\epsilon)\frac{2^{k+1}}{ek}$ for a
fixed $\epsilon>0$, then a more careful analysis shows that the height of the
tree created becomes $O_{\epsilon}(k)$. This means that the size of the tree and the
corresponding formula is bounded by a {\em polynomial in $d$}.
Here we just sketch an informal argument 
for this statement, using the continuous approach of
Section~\ref{sec:Informal}. 
We briefly recall and follow through the idea of
Subsection~\ref{subsection:cont}, keeping track of the height of the
constructed tree and neglecting the same
small error terms neglected there. Our choice of $d$ translates to the
setting of the parameter $T=\frac{e}{1+\epsilon}$ there. We build a
$(k,d)$-tree starting from a full binary tree of depth $k$. (In fact
the true depth is $k-o(k)$, but we neglect this difference too.) We
associate vectors of length $k+1$ (the leaf-vectors) with the
vertices. In the process we ignore the leaves whose leaf-vectors
satisfy the conditions of Lemma~\ref{lemma:payofflemma}, as they are
constructible. At any given step of the construction all other leaves
share the same leaf-vector and their distance from the root is the
depth of the current tree. At the next step we attach new trees to the
non-ignored leaves and thereby increase the depth of the current tree. 
We approximate the development of the leaf-vectors of the leaves
with a two variable function $F(t,x)$, where $t$ refers to ``time''
and, in fact, corresponds to the depth of the current tree (not counting
the depth of the initial full binary tree). More precisely, the $j$th
term of the leaf-vector of vertices in depth $m>k$ is approximated by
$F\left(\frac{m}{k}-1,\frac{j}{k}\right)\frac{d}{2^{k+1-j}}$. In Subsection~\ref{subsection:cont}  we develop a
differential equation for $F$ and by studying it, we show that
$\int_0^1F(t,x)\enspace dx>T$ is satisfied for some value of $t$
depending only on $T<e$ (and, through $T$, depending on $\epsilon$).  This means that the leaf-vectors at distance $tk+k$ from the root satisfy the conditions of Lemma~\ref{lemma:payofflemma}. We stop this process at this depth and obtain a binary tree with {\em all} its leaves satisfying the requirements of the lemma. This tree is still incomplete, but the construction of the $(k,d)$-tree can be completed with attaching suitable binary trees to all the leaves of this incomplete tree. The existence of suitable trees is guaranteed by Lemma~\ref{lemma:payofflemma} and the simple proof of the lemma shows that these trees can be chosen to have depth at most $k$.
This concludes the argument, which then can be made formal using a discretization
analogous to the one in Section~\ref{sec:FormalProof}.

Let us define $f_1(k,d)$ for $d>f_1(k)$ to be the minimal size of a
$(k,d)$-CNF in MU(1), and let $f_{\tree}(k,d)$ stand for the minimal size of a
$(k,d)$-tree, assuming $d\ge\f(k)$. While $f_1(k,f_1(k)+1)$ and similarly
$f_{\tree}(k,f_{\tree}(k)))$ are probably doubly exponential in $k$, the
above argument shows that for slightly larger
values of $d=(1+\epsilon)f(k)$ the values $\f(k,d)$ and hence also $f_1(k,d)$ 
are polynomial in $d$ (and thus simply exponential in $k$).

\subsection{Extremal Values and Algorithms}

Finally, we mention the question whether $\f(k)=f_1(k)+1$ for all
$k$. (The $+1$ term comes from us defining these functions
inconsistently. While $f(k)$ and $f_1(k)$ is traditionally defined as
the {\em largest} value $d$ with all $(k,d)$-CNF satisfiable or no $(k,d)$-CNF
in MU(1) exist, respectively, it was more
convenient for us to define $\f(k)$ as the {\em smallest} value $d$ for which
$(k,d)$-trees exist.) This question asks if one necessarily loses (in the
maximal number of appearances of a variable) by selecting
the $k$ vertices farthest from the root when making an MU(1) $k$-CNF formula
from a binary tree. As mentioned above,
$f(k)=f_1(k)$ is also open, but $\f(k)=f_1(k)+1$ seems to be a
simpler question as both functions are computable. Computing their values up
to $k=8$ we found these values agreed. To gain more insight we computed the
corresponding size functions too and found that $\f(k,d)=f_1(k,d)$ for
$k\le7$ and all $d>f_1(k)$ with just a single exception. We have $f_1(7)=17$
and $f_1(7,18)=10,197,246,480,846<\f(7,18)=10,262,519,933,858$.
We also found $f_1(8,d)=\f(8,d)$ for $f_1(8)=29<d\le33$.
Is $f_1(7,18)<\f(7,18)$ the exception or will it turn into the rule for
larger values? Does it indicate that $\f(k)$ and $f_1(k)+1$ will eventually
diverge?

A related algorithmic question is whether the somewhat simpler structure of
$(k,d)$-trees can be used to find an algorithm computing $\f(k)$ substantially
faster than the algorithm of Hoory and Szeider \cite{HSCompMU} for computing
$f_1(k)$. Such an algorithm would give useful estimates for $f_1(k)$ and also
$f(k)$. At present we use a similar (and similarly slow) algorithm for
either function.

\subsection{Pairing strategies and the Neighborhood Conjecture}\label{newest}
Determining 
whether the Local Lemma $2$-coloring bounds in \eqref{lll1} and \eqref{lll2} 
have any game-theoretic generalizations is still open. 
To formalize this problem we introduce
$$D(k) := \min \{ d:  \mbox{ $\exists$ $k$-uniform Maker's win ${\cal F}$ with $\Delta ({\mathcal F})=d$} \}$$
as the critical maximum degree for the Neighborhood Conjecture.
By definition, every $k$-uniform ${\cal F}$ with $\Delta ({\cal F}) \leq D(k)-1$ 
is Breaker's win and hence it is $2$-colorable by (\ref{eq:breaker-2-color}). 
 
In \cite{Beckbook} various weaker versions of the original
Neighborhood Conjecture are stated, maybe the most fundamental one is
whether $D(k)$ grows exponentially:
\begin{problem}
Does there exist an $\epsilon>0$ such that
$D(k)>(1+\epsilon)^k$ holds for all $k$?
\end{problem}

Theorem~\ref{maker-breaker} provides an upper bound for $D(k)$. For the proof of this theorem we
constructed a $k$-uniform hypergraph with bounded maximum degree and gave Maker a
winning strategy based on a pairing of the vertices. 
The point we would like to make in this section is that in order to approach
the Neighborhood Conjecture any closer, the study of these kind of 
pairing strategies is  {\em insufficient} in either direction: we {\em must} develop
methods that are able to tackle non-pairing strategies.

We call a strategy for either player in the Maker-Breaker game a {\em
pairing strategy} if it is defined by a partition of a subset of the
vertex set into pairs and calls for the player to respond with
claiming the pair of the vertex last claimed by her opponent. (In case
the opponent claimed an unpaired vertex or a vertex whose pair is
already claimed, the move is arbitrary.) In addition, a pairing strategy
for Maker also specifies a starting vertex, disjoint from the pairs, to be claimed first.
Note that a player playing by a pairing strategy is guaranteed to
eventually claim at least one member of every pair.
We call a hypergraph {\em PairingBreaker's win} (respectively, {\em PairingMaker's
win}) if Breaker (respectively, Maker) has a winning pairing
strategy.

Equivalently, we call a hypergraph $(V, {\cal F})$ a {\em PairingBreaker's win}, 
if there exists a partition of a subset 
$V' = \cup  X_i \subseteq V$ into disjoint subsets $X_i$ of size two,
such that for every $F\in {\cal F}$ there exists an index $i$ with $X_i \subseteq F$. 
The hypergraph is a  PairingMaker's win, if there
exists a vertex $x_0\in V$ and a partition of a subset 
$V' = \cup  X_i \subseteq V\setminus \{ x_0\}$ into disjoint subsets $X_i$ of size two,
such that for every subset $M$, $x_0 \in M \subseteq V$, with $|M\cap X_i| \geq 1$ 
for every $i$, there exists an $F\in {\cal F}$ with $F \subseteq M$. 
To see the equivalence note that for partitions that leave out at most a single element of
$V$, the corresponding pairing strategy is winning if and only if the above 
combinatorial condition is satisfied.

Analogously for $D(k)$, we define
\begin{align*}
D_{\pairing}(k) & := \min \{ d:  \mbox{ $\exists$ $k$-uniform
 PairingMaker's win ${\cal F}$ with $\Delta ({\mathcal
   F})=d$} \},\\
D^*_{\pairing}(k) & := \min \{ d:  \mbox{$\exists$ $k$-uniform  ${\cal F}$ with $\Delta ({\mathcal
   F})=d$, not PairingBreaker's win}\}.
   \end{align*}

Requiring a player to play a pairing strategy is a restriction, so we
have:
\begin{equation}\label{dk}
D^*_{\pairing}(k)\leq D(k)\le D_{\pairing}(k).
\end{equation}

It is interesting to note that both the functions $D^*_{\pairing}$
and $D_{\pairing}$ are well understood, but $D$ is not. The
above inequalities are the only known bounds for $D(k)$. 

The value of $D^*_{\pairing}(k)$ can be easily concluded from the classic 
observation of Hales and Jewett~\cite{hales-jewett} connecting Hall's Theorem to 
pairing strategies (cf. \cite{Beckbook}).
We indicate the simple proof to be self-contained.

\begin{prop} For every $k$ we have $D^*_{\pairing}(k)=\lfloor
k/2\rfloor+1$.
\end{prop}

\begin{proof} 
For the lower bound 
let us assume first that $\cal F$ is a $k$-uniform hypergraph with
$d:=\Delta({\cal F})\le k/2$. By Hall's theorem we can select two
representatives $x_F,y_F\in F$ for every edge $F\in\cal F$ such that
all these $2|{\cal F}|$ elements are distinct. These pairs provide Breaker
with a winning pairing strategy, since he can occupy one vertex from each
element of ${\cal F}$. 

For the upper bound let us consider any
$k$-uniform $d$-regular hypergraph $(V,{\cal F})$ with 
no two hyperedges sharing more than a single vertex. 
One can, 
for example, take the vertex set $V=\{ 1, \ldots , k\}^d$ to be the $d$-dimensional grid, 
and the family $\mathcal{F}$ consisting of all its subsets on axis parallel lines.
We show $D^*_{\pairing}(k)\le\lfloor k/2\rfloor+1$ by showing that 
the hypergraph $(V,{\cal F})$ is not PairingBreaker's win if $d > k/2$. 
Indeed, if Breaker had a winning pairing strategy, then there is a family of disjoint
pairs such that all hyperedges 
$F\in\cal F$ contain one of these pairs. These pairs must be distinct for different 
hyperedges as two hyperedges of ${\cal F}$ do not share two vertices. 
But then we have $|V|/2\ge|{\cal
F}|$ and therefore $k/d=|V|/|{\cal F}|\ge2$, providing a contradiction.
\end{proof}

Our results in this paper determine $D_{\pairing}(k)$ asymptotically.
\begin{thm}For every $k\ge3$ we have
$$f(k)/2<D_{\pairing}(k)\le2f_\tree(k-2).$$
In particular, we have
$$D_{\pairing}(k) = \frac{2^k}{ek}(1+o(1)).$$
\end{thm}
\begin{proof}
The upper bound follows from part $(i)$ of Theorem~\ref{maker-breaker} 
as the winning strategy we construct there for Maker is, in fact, a
pairing strategy. Actually, the slightly stronger bound
$D_\pairing(k)\le f_\tree(k-1)$ can also be proved with the same
argument.

For the lower bound we construct an unsatisfiable $(k,2\Delta({\cal H}))$-CNF formula $F({\cal H})$ 
from a PairingMaker's win hypergraph $(V, {\cal H})$. 

First we create a new hypergraph ${\cal F}$ by taking two disjoint copies 
of ${\cal H}$ and define a partition of its vertex set into pairs by
keeping the pairs from the pairing strategy in both copies and adding
a last pair consisting of the starting vertices in the two copies.
Then, for each pair $Y_i$ in our partition $\cup Y_i \subseteq V({\cal F})$ we introduce a variable $v_i$
and label one of the vertices in $Y_i$ with $v_i$ and the other with $\bar{v}_i$. 
For the rest of the vertices, from $V({\cal F}) \setminus \cup Y_i$, we introduce 
a separate new variable for each.
We construct a clause $C(F)$ from each hyperedge 
$F\in {\cal F}$ by taking the disjunction of the labels of the vertices.
Finally $F({\cal H})$ is the conjunction of the clauses $C(F)$ for $F\in\cal F$. 
Note that each variable
appears at most $2\Delta ({\cal H})$ times in $F({\cal H})$. 

Consider now any assignment
$\alpha$ to this formula. We want to prove that $\alpha$
is not a satisfying assignment.
Let us consider the set $S$ of vertices of $\cal F$ whose label in $\alpha$ evaluates
to {\tt false}. Then $S$ contains the starting vertex of exactly one of the copies ${\cal H}_1$ of $\cal H$, as well as one vertex from each pair in this copy. Since the pairing
defines a winning strategy for Maker in ${\cal H}_1$ there exists a set $F\in {\cal H}_1$ such that
$F\subseteq S$. Each literal of the corresponding clause $C(F)$ evaluates to 
${\tt false}$, so $\alpha$ is not satisfying.

Finally, the asymptotic statement follows from our lower and upper
bounds, and Theorems~\ref{main} and \ref{kd-trees}, respectively.
\end{proof} 

The gap between the lower and upper bounds for $D(k)$ in
Inequality~(\ref{dk}) is huge. To make it any smaller, by the above, one {\em must}
consider non-pairing strategies. Even the following problem is open:

\begin{problem} \label{pr:Hall}
Does $D(k)=D^*_\pairing(k)=\lfloor k/2\rfloor+1$ hold for every $k$?
\end{problem}

The equality above holds for $k\le4$ as the hypergraph dual of the Petersen 
graph is a Maker's win 3-uniform hypergraph
of maximum degree 2 and Knox~\cite{Knox} constructed 4-uniform Maker's
win hypergraphs with maximum degree $3$.
However Problem~\ref{pr:Hall} is open for every
$k>4$. It seems likely that resolving it in either direction 
will require a new idea and hence might lead
to a more significant progress on the Neighborhood Conjecture itself.

A more modest problem is to try to separate
$D(k)$ from $D_{\pairing}(k)$. We know that $D(2)=2=D_{\pairing}(k)$, 
as shown by a triangle. What happens for $k=3$?
We of course know that the 
hypergraph dual of the Petersen graph is a $3$-uniform Maker's win hypergraph
with maximum degree $2$, but does there exists a $3$-uniform 
hypergraph with maximum degree $2$ which is a PairingMaker's win? 

\begin{problem} 
Prove that for all large enough $k$ we have $D(k) < D_{\pairing}(k)$.
\end{problem}

{\bf Acknowledgment:}  We would like to thank the referees for their many insightful
comments.

\end{document}